\documentclass[11pt]{amsart}

\newcommand{\BZ}{\mathbb{Z}}
\newcommand{\BN}{\mathbb{N}}
\newcommand{\BR}{\mathbb{R}}
\newcommand{\BT}{\mathbb{T}}

\newcommand{\BC}{\mathbb{C}}

\newcommand{\BP}{\mathbb{P}}

\newcommand{\caB}{\mathcal{B}}

\newcommand{\caD}{\mathcal{D}}
\newcommand{\caF}{\mathcal{F}}

\newcommand{\caG}{\mathcal{G}}

\newcommand{\caH}{\mathcal{H}}

\newcommand{\gC}{\Gamma}

\newcommand{\gL}{\Lambda}
\newcommand{\gF}{\Phi}
\newcommand{\ga}{\alpha}

\newcommand{\gd}{\delta}

\newcommand{\gc}{\gamma}
\newcommand{\gf}{\varphi}
\newcommand{\gl}{\lambda}
\newcommand{\gs}{\sigma}
\newcommand{\gr}{\rho}
\newcommand{\gx}{\chi}
\newcommand{\gp}{\psi}

\newcommand{\gz}{\zeta}

\newcommand{\bt}{\bar{\triangle}}

\newcommand{\ti}[1]{\tilde{#1}}
\newcommand{\p}{\prod}

\newcommand{\te}{\text}

\newcommand{\ra}{\rightarrow}

\newcommand{\lra}{\longrightarrow}

\newcommand{\ve}[1]{\stackrel{\ra}{#1}}
\newcommand{\converge}[1]{\stackrel{#1}{\lra}}

\newcommand{\ip}[2]{\left<{#1},{#2}\right>}
\newcommand{\Lim}[1]{\lim_{{#1} \ra \infty}}

\newcommand{\Sum}[2]{\sum_{{#1}=1}^{{#2}}}
\newcommand{\Avr}[2]{\frac{1}{{#1}}\sum_{{#2}=1}^{{#1}}}

\swapnumbers
\theoremstyle{plain}
\newtheorem{lma}{Lemma}[section]
\newtheorem{thm}[lma]{Theorem}
\newtheorem{pro}[lma]{Proposition}

\newtheorem{cor}[lma]{Corollary}

\theoremstyle{definition}
\newtheorem{dfn}[lma]{Definition}
\newtheorem{rmr}[lma]{Remark}
\newtheorem{ntt}[lma]{Notation}

\newtheorem{dsc}[lma]{}

\newcommand{\drorisgreat}[1]{}

\begin{document}
\title[Universal Characteristic Factors]{Universal Characteristic Factors
and Furstenberg Averages} 
\author{T. Ziegler }

\email{ tamar@math.ohio-state.edu}
\address{Department of Mathematics, The Ohio State University, 
         Columbus, Ohio}

\subjclass[2000]{37Axx}
\maketitle

\begin{abstract}
Let $X=(X^0,\caB,\mu,T)$ be an ergodic probability 
measure preserving system. For a natural
number $k$ we consider the averages
\begin{equation*}\tag{*} \Avr{N}{n} \p_{j=1}^k f_j(T^{a_jn}x)
\end{equation*}
where $f_j \in L^{\infty}(\mu)$, and $a_j$ are integers. 
A factor of $X$ 
is characteristic for averaging schemes
of length $k$ (or $k$-characteristic) if for any non zero distinct 
integers $a_1,\ldots,a_k$,
the limiting $L^2(\mu)$ behavior of the averages in (*) is 
unaltered if we first project the
functions $f_j$ onto the factor. A factor of $X$ is a 
{\em $k$-universal characteristic factor
($k$-u.c.f)} if it is a $k$-characteristic factor,
and a factor of any $k$-characteristic factor.
We show that there exists a unique $k$-u.c.f, 
and it has a structure of a $(k-1)$-step nilsystem, 
more specifically an inverse limit of $(k-1)$-step nilflows. Using this
we show that the averages in (*) converge in $L^2(\mu)$. 
This provides an alternative proof to the one given in 
Host and Kra \cite{HK02c}.
\end{abstract}

\section{Introduction}

Averages of the form
\begin{equation}\label{avr}
 \Avr{N}{n} \p_{j=1}^k f(T^{jn}x)
\end{equation}
were first introduced by Furstenberg \cite{Fu77}
in his ergodic theoretic proof of Szemer\'edi's theorem on arithmetic progressions
in sets of positive density in $\BZ$. Furstenberg proved the following theorem:
\begin{thm}[Furstenberg]
Let $X=(X^0,\caB,\mu,T)$ be a measure preserving system (m.p.s.). Let $A$ be a 
set of positive measure, and $f=1_A$. Then
\[
 \liminf_{N \ra \infty} \Avr{N}{n} \int \p_{j=0}^k f(T^{jn}x)d\mu >0.
\]
\end{thm}
The theorem above ensures that there exists an integer $n$, such that the 
points $x,T^nx,\ldots,T^{kn}x$ are in $A$, and corresponds to the 
existence of an arithmetic progression of length $k+1$  in sets of positive 
density in $\BZ$. 

The $L^2$ limiting behavior of the averages in  equation (\ref{avr}) 
is related to a natural series of factors of the measure preserving system 
$X$. The factor 
corresponding to arithmetic progression of length $3$ (the case where $k$ is
$2$) - the Kronecker factor - was described in  \cite{Fu77}. 
The factor corresponding
 to arithmetic progression of length $4$, an inverse limit of $2$-step 
nilflows, was studied by  Conze - Lesigne  
\cite{CL84},\cite{CL87},\cite{CL88}, Furstenberg - Weiss \cite{FuW96}, 
and Host - Kra  \cite{HK01},\cite{HK02a}, and hinted to the nature of the 
complete series. The complete series of factors was discovered by 
Host-Kra \cite{HK02c}, and independently, though somewhat later, by the 
author. We give an equivalent definition of the series of factors
described in \cite{HK02c}, and provide a different  
construction for these factors. Although there are some similarities
between the constructions (for example, both start out with the Furstenberg 
structure
theorem \cite{Fu77}), the definition of the factors and the bulk of the construction are 
significantly different. 
In particular, we do not use the Gowers uniformity norms 
\cite{G01}, which are fundamental in the approach of Host-Kra to this problem
(and in the works of Gowers \cite{G01}, and Green-Tao \cite{GT04} on 
problems of a similar nature). 
The averages studied in the paper are of a special kind, but the
techniques developed can be used in analyzing other multiple
averages (e.g. averages along a polynomial sequence).

Let $X=(X^0,\caB_{X},\mu_{X},T_{X})$ be a probability 
measure preserving system;
i.e. $(X^0,\caB_{X},\mu_{X})$ is a measure space, 
and $T_X$ is a measure preserving transformation. When there is 
no confusion we will omit the subscript $X$.
We write $Tf$ for the function $Tf(x)=f(Tx)$.
By ergodic decomposition it will suffice to study the limit of 
(\ref{avr}) with the additional hypothesis of ergodicity. The nature of 
the limit will depend on mixing properties of the system.
The maximal degree of mixing relevant is  {\em weak mixing};
indeed in this case Furstenberg has shown in \cite{Fu77}:
\begin{thm}[Furstenberg]\label{WM}
If $X$ is weak mixing then
\[
  \Avr{N}{n}  \p_{j=1}^k f_j(T^{jn}x) 
   \converge{L^2(X)}  \p_{j=1}^k \int f_j(x) d\mu.
\] 
\end{thm}

For a general ergodic system $X$ the averages in 
equation (\ref{avr})
need not converge to a constant 
function. Indeed, if the system ${X}$  is not weakly mixing
there exists a nontrivial eigenfunction $\gp$. If $T\gp(x)=\gl \gp(x)$ then
\[
  T^n \gp^2(x)T^{2n}\gp^{-1}(x)=\gp(x)
\]
for all $n$, thus 
\[
 \Avr{N}{n} T^n \gp^2(x)T^{2n}\gp^{-1}(x)=\gp(x).
\]
By the above equation, the set of limiting functions 
contains the algebra spanned by  eigenfunctions 
- the {\em Kronecker algebra}. The Kronecker algebra 
determines the  `Kronecker factor' ${Z}$ 
where $Z^0$ is a compact Abelian group, $\caB_Z$ the (completed) Borel algebra,
$\mu_Z$ the Haar measure, and the action of $T_Z$ is given by translation
by an element $\ga \in Z^0$, i.e $T_Zz=z+\ga$.
Let $\pi:X \ra Z$ be the factor map. 
It is not surprising that an
Abelian group factor should come up when studying the relations between
$x,T^nx,T^{2n}x$, as the projections of these points on the Abelian group 
factor
$\pi(x),\pi(x)+n\ga,\pi(x)+2n\ga$ form an arithmetic progression:
$\pi(x)=2(\pi(x)+n\ga)-(\pi(x)+2n\ga)$. 
It turns out that this `constraint' imposed by 
the Kronecker factor is the only `constraint' on the triple 
$x,T^nx,T^{2n}x$, 
and in a manner to  be made precise, the Kronecker factor is 
`characteristic' for the limit of the averages 
$\Avr{N}{n} f(T^nx)g(T^{2n}x)$.

Let ${X}$ 
be a measure preserving system (m.p.s).
Let ${Y}$ 
be a homomorphic image; i.e,  we have a map  $\pi:X^0 \ra Y^0$
with $\pi^{-1}\caB_Y \subset \caB_X$, $\pi \mu_X =\mu_Y$ and $T_Y\pi=\pi T_X$. 
Then ${Y}$ is  a {\em factor} of ${X}$, ${X}$
is an {\em extension} of ${Y}$, and abusing the notation we write
$\pi:X \ra Y$ for the factor map. A factor of $X$ is determined by 
a $T_X$ invariant subalgebra of $L^{\infty}(\mu)$.
The map $\pi$ induces two natural maps
$\pi^*: L^2(\mu_{Y}) \ra  L^2(\mu_{X})$ given by $\pi^*f=f\circ \pi$, and 
$\pi_*: L^2(\mu_{X}) \ra  L^2(\mu_{Y})$ given by $\pi_*f=E(f|\caB_{Y})$
(the orthogonal projection of $f$ on $\pi^*L^2(\mu_Y)$). 
We fix an ergodic m.p.s $X$.
The notion of `characteristic factors' was first introduced in
a paper by Furstenberg and Weiss \cite{FuW96}.
\begin{dfn} Let ${Y}$ be a factor of $X$. Let 
$k$ be a natural number, $(a_1,\ldots,a_k)$ be distinct non-zero integers.
The system ${Y}$ is {\em characteristic for 
$(a_1,\ldots,a_k)$} if for any $f_1,\ldots,f_k \in L^{\infty}(\mu_X)$, 
\[
  \Avr{N}{n} \p_{j=1}^k  T_X^{a_jn}f_j 
 - \pi^*\Avr{N}{n} \p_{j=1}^k T_Y^{a_jn}\pi_{*}f_j
  \converge{L^2(\mu_X)} 0.
\] 
The system ${Y}$ is a
{\em $k$-characteristic} factor of $X$ if it is characteristic for any 
$k$-tuple of 
distinct non-zero integers.
\end{dfn}
It is in this sense that the Kronecker factor is characteristic for 
calculating the limit of the averages 
$\Avr{N}{n} f(T^nx)g(T^{2n}x)$.
We now define a universal characteristic factor:
\begin{dfn} Let ${Y}$ be a factor of $X$.
The system ${Y}$ is a
{\em $k$-universal characteristic factor (u.c.f)} of $X$ if it is a 
$k$-characteristic 
factor of $X$ , and a 
factor of any other $k$-characteristic factor of $X$.
\end{dfn}

For the averages $\Avr{N}{n} f_1(T^nx)f_2(T^{2n}x)f_3(T^{3n}x)$, 
the Kronecker factor does not suffice. 
Let $\gf$ be a {\em second order eigenfunction}, i.e.,
$T\gf=\gp \gf$, and  $T\gp=\gl\gp$. Then
one can check that 
\[
  T^n \gf^3(x)T^{2n}\gf^{-3}(x)T^{3n}\gf(x)=\gf(x)
\]
for all $n$; thus 
\begin{equation}\label{second_order}
 \Avr{N}{n} T^n \gf^3(x)T^{2n}\gf^{-3}(x)T^{3n}\gf=\gf(x).
\end{equation}
Let $Y$ be a factor of $X$ that is characteristic for 
$(1,2,3)$.
Equation (\ref{second_order}) implies 
that the algebra generated by all second order eigenfunctions is a subset
of  $L^{\infty}(\mu_{Y})$. 
It is natural to conjecture that the algebra generated by
the second order eigenfunctions determines a factor that is 
characteristic for $(1,2,3)$. 
Furstenberg and Weiss presented the following  counter example. Let
\[ 
   X=
   \left( \begin{smallmatrix}
   1 & \mathbb{R} & \mathbb{R} \\
     &          1 & \mathbb{R} \\
     &            &          1 \\
   \end{smallmatrix}\right)\Big/
   \left( \begin{smallmatrix}
   1 & \mathbb{Z} & \mathbb{Z} \\
     &          1 & \mathbb{Z} \\
     &            &          1 \\
   \end{smallmatrix}\right)
   =N/\Gamma.
\]
Consider the system ${X}$
where $X^0=N/\Gamma$, 
$\caB_X$ the (completed) Borel algebra, $\mu_X$ the unique measure 
invariant under translations by any element of the  group $N$, and 
$T_X$ is given by 
$T_Xg\gC =ag\gC$ for some $a \in N$ acting ergodically.
This system has no second order eigenfunctions, but there are relations
between $g\gC, a^ng\gC, a^{2n}g\gC,a^{3n}g\gC$ not coming from
Kronecker factor: In $N/\Gamma$, $g\gC$ is determined by 
$a^ng\gC, a^{2n}g\gC,a^{3n}g\gC$.

 This  system can be viewed as a circle extension of the Kronecker 
factor which is a $2$ dimensional torus, and the action of $T_X$ on 
$\BT^2 \times S^1$ is given by $T_X(z,\gz)=(z+\ga,\gs(z)\gz)$ 
(the function $\gs(z)$ is called the {\em extension cocycle}). 
The projection of the points $g\gC, a^ng\gC, a^{2n}g\gC,a^{3n}g\gC$
on the Kronecker factor will form an arithmetic progression, but as 
 $g\gC$ is a function of $a^ng\gC, a^{2n}g\gC,a^{3n}g\gC$ the points
$g\gC, a^ng\gC, a^{2n}g\gC,a^{3n}g\gC$ will not be independent on the fibers 
over the Kronecker factor. This  fact translates to 
a restriction on the extension cocycle $\gs(z)$ known as the Conze-Lesigne 
equation.
(this equation is analyzed in  \cite{CL84}, \cite{CL87}, \cite{CL88}, 
\cite{Le84}, \cite{Le87}, \cite{Le93}, \cite{FuW96}, \cite{HK01} \cite{HK02a},
\cite{Me90},\cite{R93}). 
In particular any factor that is characteristic for 
$(1,2,3)$ will contain functions other than
first and  second order eigenfunctions.

In general, if $N$ is a $k$-step nilpotent group ($N_{k+1}=1$), 
$\Gamma < N$ then $x \in N/\gC$ is determined by 
$a^nx,a^{2n}x, \ldots,a^{(k+1)n}x$.

It is natural to ask whether these are the only constraints, i.e.
do {\em all} the constraints on the points $x,T^nx, \ldots, T^{(k+1)n}x$
come from a $k$-step nilpotent factor?

\begin{dfn}
A {\em nilsystem} consists
of a space $X$ on which a nilpotent group $N$ acts transitively preserving
a measure $\mu_X$, and a transformation $T_X$ which acts by translation by a 
group element $a$: $T_Xx=ax$. A special case is when $N$ is a $k$-step 
nilpotent Lie group, $\gC$ a cocompact lattice, $X=N/\gC$ (a nilmanifold),  
and $\mu_X$ the unique measure invariant under translation by elements of $N$.
We call this a $k$-step {\em nilflow}.
A {\em $k$-step pro-nilflow} is an 
inverse limit of $k$-step  nilflows.
\end{dfn}

We prove the following theorems:
\begin{thm}\label{thm:universal} Let ${X}$ be an ergodic measure preserving system. Then there exists a unique $k$-universal characteristic factor of $X$.  
If $\pi:{X} \ra {Y}$ is a factor 
map,
and $ {W}({X}), {W}({Y})$ are $k$-universal 
characteristic factors of  
${X},{Y}$
respectfully, then $\pi$ induces a map between ${W}({X})$ and ${W}({Y})$. 
\end{thm}
If we denote by ${Y}_k({X})$ the $k$-u.c.f of $X$, then one 
obtains 
an inverse series of factors $\ldots \ra {Y}_k({X}) \ra {Y}_{k-1}({X}) \ra 
\ldots  \ra {Y}_1({X})$.
\begin{thm}\label{main} Let ${X}$ be an ergodic measure preserving system, and
let ${Y}_k({X})$ be the $k$-universal characteristic factor of $X$. 
Then ${Y}_k({X})$ has a structure of a 
$(k-1)$-step nilsystem, more specifically a $(k-1)$-step pro-nilflow.
\end{thm}

The proof of Theorem \ref{main} is by induction: 
assuming that the $k$-u.c.f is a 
$(k-1)$-step nilsystem, one reduces the problem of determining the  
$(k+1)$-u.c.f
to the case where the system ${X}$ is a circle extension of a $(k-1)$- 
step nilsystem. 
If the points $x,T^nx,\ldots T^{(k+1)n}x$ are independent
on the fibers over the $k$-step nilsystem then the $k$-step nilsystem
would suffice, i.e would be $k$-characteristic.
Otherwise one would get a restriction on the extension cocycle. 
The main difficulty is using this restriction to construct a nilpotent
group acting transitively on ${X}$.

As a corollary we get a theorem proved recently by Host and Kra \cite{HK02c} 
\begin{cor} Let $X$ be a m.p.s. Let $k$ be a natural number,
$a_1,\ldots a_k \in \BZ$, and $f_1 \ldots f_k \in L^{\infty}(\mu_X)$
then the averages 
\begin{equation}\label{eq:conv} \Avr{N}{n} \p_{j=1}^k f_j(T^{a_jn}x)
\end{equation}
converge in $L^2(\mu_X)$.
\end{cor}
By Theorem \ref{main}, in order to have $L^2(\mu_X)$ convergence of the averages in 
(\ref{eq:conv}), it is enough to prove an 
$L^2(\mu_X)$ convergence theorem for $k$-step pro-nilflows.
 Convergence in $L^2$
for pro-nilflows follows from convergence for nilflows. For nilflows
one has a.e. convergence (\cite{P69}, \cite{Sh96}, \cite{L02}). 
An explicit description of the limit is given in \cite{Le89} for the case
$k=3$, and in general in \cite{Z02a}.

{\bf Acknowledgment}
Most of the ideas in this work appear in the authors PhD thesis. I
would like to thank my adviser Prof. Hillel Furstenberg for introducing
me to ergodic theory, specifically to questions involving non-conventional
ergodic averages, and for many fruitful discussions. I would also like
to thank Benji Weiss, Shahar Mozes, Vitaly Bergelson for
enlightening conversations and valuable remarks. I owe special thanks
to Sasha Leibman for pointing out many inaccuracies in the early version
of this paper.

\section{ Universal characteristic factors}

We start by proving Theorem \ref{thm:universal}.
\begin{lma}
Let $X$ be an ergodic m.p.s.
Let ${Y}_1,{Y}_2$ be $k$-characteristic factors of $X$.
Then there exists a $k$-characteristic factor of $X$, which
is a factor of both ${Y}_1,{Y}_2$.
\end{lma}
\begin{proof}
Denote $P,Q$ the orthogonal projections onto 
$L^2(\mu_{{Y}_1}),L^2(\mu_{{Y}_2})$ 
(seen as subspaces of $L^2(\mu_{X})$) respectfully, and by $\pi_i:X^0 \ra Y^0_i$
for $i=1,2$ the factor maps.
Then $P^2=P^*=P$ (same for $Q$).
We show that  $(PQP)^n$ strongly converges to a self adjoint operator 
projection $W$: $P$ is a projection thus $P \le I$. 
\[
  \ip{(PQP)^2x}{x}=\ip{PQPx}{QPx}\le \ip{QPx}{QPx}=\ip{PQPx}{x},
\]  
Inductively, the sequence $(PQP)^n$ is a decreasing sequence
of operators, thus $\ip{(PQP)^nx}{x}$ converges for all $x$. The sequence
$(PQP)^nx$ is a Cauchy sequence as
\[\begin{aligned}
  \|(PQP)^nx-(PQP)^mx\|^2=&\ip{(PQP)^{2n}x}{x}+\ip{(PQP)^{2m}x}{x}\\
                          &-2\ip{(PQP)^{(n+m)}x}{x} \ra 0.
  \end{aligned}
\]
Let $W=\lim_{n \ra \infty} (PQP)^n$, then $W^2=W=W^*$. If $Wx=x$
then $Px=PWx=Wx=x$, and
\[
 PQx=PQPx=PQPWx=Wx=x \Rightarrow \|Qx\|=\|x\| \Rightarrow Qx=x. 
\]
It follows that $W(L^2(X^0,\caB_X,\mu_X))= L^2(X^0,\caD,\mu_X)$ for 
$\caD=\pi_1^{-1}(\caB_{Y_1}) \cap \pi_2^{-1}(\caB_{Y_2})$.
We show that ${W}(L^2(\mu_{X}))$ is a $k$-characteristic factor of $X$. 
For all $m$:
\[
 \begin{aligned}
  \lim_{N  \ra \infty} &\Avr{N}{n} T^{a_1n}f_1\ldots T^{a_kn}f_k \\
 =&\lim_{N \ra \infty}
   \Avr{N}{n} T^{a_1n}((PQP)^mf_1)\ldots T^{a_kn}((PQP)^mf_k)
 \end{aligned}
\]
\end{proof}

\begin{cor} Let $X$ be a m.p.s. 
There exists a unique $k$-universal characteristic factor of $X$.
\end{cor}

\begin{proof}
By Zorn's lemma.
\end{proof}

The advantage of looking at all $k$-tuples (rather than focusing on a 
specific one )
is that $k$-u.c.f are  natural in the sense that
any morphism of  measure preserving systems  induces a morphism 
between their $k$-universal characteristic factors - 
as will be shown in corollary \ref{natural}.
(This may also be true for characteristic factors of a specific scheme).

\begin{lma}\label{lma:limit_char} Let $V$ be the algebra generated by partial limits
of the sequences $\{\Avr{N}{n} T^{a_1n}f_1\ldots T^{a_kn}f_k \}$, 
where $f_i \in L^{\infty}(\mu_X)$, and $a_0,\ldots,a_k \in \BZ$ are distinct 
non zero integers.
Then $V$ determines the $k$-universal characteristic factor of $X$. 
\end{lma}

\begin{proof} Let $W(X)$ be the $k$-universal characteristic factor. 
Obviously $V \subset L^{\infty}(\mu_{W(X)})$.
We show that the factor determined by 
$V$ is a $k$-characteristic factor of $X$.
Let $g \perp V$, then for any $f_1$ 
\[ \begin{aligned}
 \ip{g}{\Avr{N}{n} T^{a_1n}f_1\ldots T^{a_kn}f_k}
 =&\Avr{N}{n} \int g T^{a_1n}f_1 T^{a_2n}f_2 \ldots T^{a_kn}f_k d\mu\\
 =&\Avr{N}{n} \int f_1 T^{-a_1n}g\ldots T^{(a_k-a_1)n}f_k d\mu \\
 &\converge{N \ra \infty} 0.
 \end{aligned}
\]
\end{proof}

\begin{cor}\label{natural}
 If $\pi:{X} \ra {Y}$ is a factor map,
and ${W}({X}),{W}({Y})$ are $k$-u.c.f. 
for ${X},{Y}$
respectfully, then $\pi$ induces a map between ${W}({X})$ and ${W}({Y})$. 
\end{cor}

\begin{dsc}{\bf Universal characteristic factors for $k=1,2$}. 
Let $X$ be an ergodic m.p.s. 
If the system ${X}$ is totally ergodic then by 
von Neumann's theorem the trivial
system (a point) would be the $1$-u.c.f of $X$.
Otherwise one needs to take into account the algebra generated by  
functions that are invariant under 
$T_X^m$ for some $m$. The factor may then be represented as a finite cyclic 
group with addition of one, or a pro-cyclic group which is an inverse limit of cyclic groups.  The $2$-universal characteristic factor of $X$ 
coincides with the first block in Furstenberg's structure 
theorem (see \cite{Fu77}) and is referred to as {\em the Kronecker factor}.
The system  ${Z}$ 
is a {\em Kronecker system} (or an 
\emph{almost periodic system}) if $Z^0$ is 
a compact abelian group (a $1$-step nilpotent group), 
$\caB_Z$ is the (completed) Borel algebra, $\mu_Z$ is the Haar 
measure, and the action of $T_Z$ is given by $T_Zz=z+\ga$ for some 
$\ga \in Z^0$. 
The Kronecker factor is the maximal almost periodic factor. 
Equivalently, ${Z}$ is the Kronecker factor of ${X}$  
if the eigenfunctions of $T_X$ span $L^2(\mu_{Z})$ (thought of as a subspace
of $L^2(\mu_{X})$). 
\begin{rmr} If the system ${X}$ is weak mixing, i.e. has no 
non trivial eigenfunctions, then the Kronecker factor is trivial
(and $Y_k(X)$ is trivial for all $k$).
\end{rmr}
\end{dsc}

\begin{dsc}{\bf Isometric extensions}\label{isometric_extensions}
The notion of characteristic factors was motivated by Furstenberg's structure
theorem \cite{Fu77}.  Furstenberg's idea was to relativize the notion of weak mixing to 
a {\em weak mixing extension} and to define the complementary notion of a 
{\em compact extension}  (or {\em isometric extension}). 
Let ${X}$ 
be an ergodic m.p.s., and let 
${Y}$ 
be a factor.
Consider the ring $L^{\infty}(\mu_{Y})$ as a subring of functions on $X$. A
subspace $V \subset L^2(\mu_{X})$ is a {\em finite rank module} over
$L^{\infty}(\mu_{Y})$ if there exist finitely many functions 
$\gf_1,\ldots,\gf_k$, such that any function $f \in V$ can be expressed
as  $f=\sum_{i=1}^k a_i(y)\gf_i(x)$. 
We say that ${X}$ is an {\em isometric extension} of ${Y}$
if $L^2(\mu_{X})$ is spanned by finite rank $T_X$ invariant modules over 
$L^{\infty}(\mu_{Y})$.
It can be shown that in this case  ${X}$ is  isomorphic to a 
{\em skew product} ${X}'$ where $X'^{0}=Y^0 \times M$, where  $M=G/H$ is a 
homogeneous compact metric space, $\mu_{X'}=\mu_Y \times m_M$, where
$m_M$ is the unique probability 
measure invariant under the transitive group of isometries $G$, 
and the action of $T_{X'}$ is given by 
$T_{X'}(y,m)= (T_Yy, \gr(y)m)$, where $\gr: Y^0 \ra G$. We denote
$T_{X'}$ by $T_{Y,\gr}$, or if there is no confusion, just $T_{\gr}$. 
For example, a Kronecker system is an isometric extension of a point.
Define $\gr^{(n)}:Y^0 \ra G$ by $T_{\gr}^n(y,m) =(T^n y, \gr^{(n)}(y)m)$; 
then $\gr^{(n)}$ satisfies
a $1$-cocycle equation for the action of $\BZ$ on functions from $Y$ to $G$
\[
  \gr^{(n+m)}(y)=\gr^{(n)}(T^my)\gr^{(m)}(y).
\]
Since $\gr^{(n)}(y)$ is determined by $\gr^{(1)}(y)$ we shall focus on 
$\gr(y)=\gr^{(1)}(y)$
and refer to it as the {\em extension cocycle} (or just {\em cocycle}).
Abusing the notation we denote the system ${X}'$ by 
${Y} \times_{\gr} G/H$.
For more details see \cite{Fu77}, or \cite{Zi76}
\end{dsc}

\begin{dsc}
Let $X_1$, $X_2$ be m.p.s. and let $Y$ be a common factor with
$\pi_i:X^0_i \ra Y^0$ for $i=1,2$ the factor maps. Let 
$\mu_{X_i,y}$ represent the disintegration of $\mu_{X_i}$ with respect to 
$Y$.
Denote $\mu_{X_1}\times_Y \mu_{X_2}$ the measure  defined by
\[
 \mu_{X_1}\times^{}_Y \mu_{X_2}(A)=\int \mu_{X_1,y}\times^{}_Y \mu_{X_2,y}(A) d\mu_Y
\] 
for $A \in \caB_{X_1} \times \caB_{X_2}$. 
The system 
\[ (X^0_1 \times X^0_2, \caB_{X_1} \times  \caB_{X_2}, 
   \mu_{X_1}\times^{}_Y \mu_{X_2}, T_{X_1} \times T_{X_2}  )\] 
is called the relative product of $X_1$ and $X_2$ with respect to $Y$ and is 
denoted ${X}_1 \times_Y {X}_2$. 
\end{dsc}

\begin{dsc}
Let ${X}$ be an ergodic m.p.s.,
${Y}$ a factor and $\pi:X \ra Y$ the factor map.
Consider the subspace of $L^2(\mu_{X})$ spanned by all finite 
rank $T_X$-invariant 
modules over  $\pi^*L^{\infty}(\mu_{Y})$. 
This subspace will be defined by some factor
$\hat{{Y}}$ between ${X}$ and ${Y}$. The system
$\hat{{Y}}$
is called the 
{\em maximal isometric extension of $Y$ in $X$}. For some $l \in BN$, let 
$X'=(X^0,\caB_X,\mu_X,T_X^l)$, and let $Y'=(Y^0,\caB_Y,\mu_Y,T_Y^l)$. 
Then the maximal isometric extension of $Y'$ in $X'$ is 
$\hat{Y'}=(\hat{Y}^0,\caB_{\hat{Y}},\mu_{\hat{Y}},T_{\hat{Y}}^l)$.
\end{dsc}

\begin{dsc}\label{conditional}
Let ${X}_i$, $i=1, \ldots,k$, be measure preserving systems,
and let ${Y}_i$ be corresponding factors, and $\pi_i:X^0_i \ra Y^0_i$ the 
factor maps. A measure $\nu$ on $\Pi Y^0_i$ defines a {\em joining} of the 
measures on $Y_i$ if it is invariant under 
$T_{Y_1} \times \ldots \times T_{Y_K}$ and maps onto $\nu_{Y_i}$ under 
the natural map  $\Pi Y_i \ra Y_j$. Let $\nu$ be a joining of the measures on 
$Y_i$, and let $\mu_{X_i,y_i}$ represent the  disintegration of 
$\mu_{X_i}$ with respect to $Y_i$.

Let $\mu$ be a measure on  $\Pi X^0_i$ defined by
\[
 \mu= \int \mu_{X_1,y_1}\times  \ldots \times \mu_{X_k,y_k} d\nu(y_1,\ldots,y_k).
\]
Then  $\mu$ is called the {\em conditional product measure} with 
respect to $\nu$.

The following is shown in \cite{Fu77}  Theorem $9.5$:
\begin{thm}[Furstenberg]\label{max_isometric}
Let ${X}_i,{Y}_i$, $\nu$, $\mu$  be as in \ref{conditional}. Assume each 
${X}_i$ has finitely many 
ergodic components. Let $\hat{{Y}}_i$ be the 
maximal isometric extension of ${Y}_i$ in 
${X}_i$, $\hat{\pi_i}:X_i \ra \hat{{Y}}_i$ the projection. 
Then if $F \in L^2(\mu)$ is invariant under 
$T_{X_1} \times \ldots \times T_{X_k}$ then there exists a function 
$\gF \in L^2(\mu)$, so that 
\[
F(x_1, \ldots,x_k)=\gF(\hat{\pi}_1(x_1),\ldots,\hat{\pi}_k(x_k)).
\]
\end{thm}
\end{dsc}

\begin{dsc}{\bf Group extensions}
A special case of an isometric extension ${X} \ra {Y}$ is when 
the homogeneous space $M$ from \ref{isometric_extensions} is equal to $G$, i.e 
$X=Y\times_{\gr} G$
where $G$ is a compact group. In this case we say that ${X}$ is a 
{\em group extension} of ${Y}$.
\end{dsc} 

\begin{lma}\label{isometric_to_group}
Suppose ${X}$ is an ergodic isometric extension of ${Y}$
so that we can express $X=Y \times_{\gr} G/H$. Using the function $\gr$, we 
can define a group extension ${Y} \times_{\gr} G$. Then $G$ and $H$
can be chosen so that the extension ${Y} \times_{\gr} G$ is an ergodic group
extension.  
\end{lma}
\begin{proof} \cite{FuW96} lemma $7.2$.
\end{proof}

\begin{lma}\label{lma:intermediate}
Let ${X}={Y} \times_{\gr} G$ be an ergodic group extension of ${Y}$, and 
let ${W}$ be an intermediate factor between ${X}$ and ${Y}$, then
${X}$ is a group extension of ${W}$.  
\end{lma}

\begin{proof} \cite{FuW96} lemma $7.3$ .
\end{proof}

\begin{dsc}\label{Mackey}
Let ${Y}$ be an ergodic m.p.s., $G$ a compact metrizable group.
Let ${Y} \times_{\gr} G$ be a group extension.
We can parameterize $Y^0 \times G$ replacing
$(y,g)$ with $F(y,g)=(y,f(y)g)$ for some measurable function $f:Y^0 \ra G$.
Let $\gr'(y):=f(Ty)\gr(y)f(y)^{-1}$, then the systems 
${Y} \times_{\gr} G$, ${Y} \times_{\gr'} G$ are isomorphic,
and $\gr,\gr'$ are called {\em equivalent cocycles} or {\em cohomologous cocycles}. If $\gr$ is equivalent
to the identity cocycle then $\gr$ is a {\em $Y$-coboundary}
(or just {\em coboundary} when there is no confusion). If  $\gr$
is equivalent to a constant cocycle then $\gr$ is a {\em $Y$-quasi-coboundary}
(or {\em quasi-coboundary}).
\end{dsc}

\begin{dsc}
If $\gr$ takes values in a closed subgroup $H$ of $G$, the
extension ${Y} \times_{\gr} G$ will not be ergodic 
(any function on $H\backslash G$ 
will be invariant). By the foregoing discussion if $\gr$ is equivalent to 
a cocycle taking values in a closed subgroup $H$, then the extension
${Y} \times_{\gr} G$ will not be ergodic. 

\begin{thm}[Mackey]\label{mackey_thm}
Let $\gr:Y^0 \ra G$ be a measurable cocycle. 
There exists  a closed subgroup $M<G$, unique up to conjugacy, so that:
\begin{enumerate}
 \item $\gr$ is equivalent to a cocycle $\gr'$ taking values in $M$.\\
       i.e. $\gr'(y)=f(Ty)\gr(y)f(y)^{-1} \in M$.
 \item Any ergodic $T_{\gr'}$-invariant measure on $Y^0 \times G$, extending $\mu_Y$,
       has the form
       $\mu_Y \times m_{M\gc}$ for some coset $M\gc$, and the ergodic $T_{\gr}$
       invariant measures are obtained by applying $F^{-1}$  
       (defined in \ref{Mackey} to the ergodic
       $T_{\gr'}$-invariant measures. 
        The group $M$ is called the {\em Mackey
       group} of the extension ${Y} \times_{\gr} G$.  
\end{enumerate}
\end{thm}
\begin{lma}\label{lma:conjugate} For $i=1,2$, let $Y_i$ be ergodic m.p.s, 
let $X_i=Y_i \times_{\gr_i} G$ be group extensions, and let $M_i$ be the
associated Mackey groups.
Let $\pi_i$ be the projection $\pi_i:X_i \ra Y_i$. Let $S:X_1 \ra X_2$ be an 
isomorphism such that $S:Y_1 \ra Y_2$ ans $S\pi_1=\pi_2S$.  
Then $M_1$ and $M_2$ are conjugate.   
\end{lma}
\begin{proof} The transformation $S$ maps the ergodic components of the 
group extension $Y_1 \times_{\gr_1} G$ onto those of $Y_2 \times_{\gr_2} G$.
Each ergodic component is determined by a right coset $M_i\gc_i$ for 
$\gc_i \in G$, thus $S$ induces a map from 
$\gf:M_1\backslash G  \ra M_2\backslash G$, that commutes with the action of 
$G$ from the right. Thus $\gf$ is a $G$-isomorphism, and therefore $M_1$ and $M_2$ are conjugate.
\end{proof}
\end{dsc}

\section{ abelian extensions}
\begin{ntt} We use additive notation for abelian groups with the 
exception of the group $S^1=\{\gz \in \BC:|\gz|=1\}$ which will play a special role in the future.
In particular, if $\gr,\gr'$ are equivalent cocycles 
(defined in the foregoing section) taking values in an abelian group $G$, then
there exists a function $f:Y^0 \ra G$ such that  
\[
  \gr(y)= f(Ty)+\gr'(y)-f(y).
\]
\end{ntt}

\begin{dsc}\label{Mackey_abelian}
Let $G$ be a compact abelian group, then ${Y} \times_{\gr} G$ is an 
{\em abelian extension}. In this case  the Mackey group 
defined in the foregoing section $M$ is unique.
Let  
\[
  M^{\perp}=\{ \gx \in \hat{G}: \gx(g)=1 \ \te{for all} \ g \in M \}
\]
be the annihilator of $M$. If $\gr$ is equivalent to a cocycle taking values
in $M$ then  $\gx \circ \gr$ is a coboundary for all
$\gx \in \hat{M}$. 
\[
  M^{\perp}=\{ \gx \in \hat{G}: \gx \circ \gr \ \te{is a coboundary} \}.
\] 
\end{dsc}

\begin{pro}\label{P:reduction_abelian}
Let ${Y} \times_{\gr} G$ be an abelian extension,
and let $M$ be the Mackey group of this 
extension. Let $f \in L^2 (\mu_{Y} \times m_G)$ be s.t. 
for all $\gx \in M^{\perp}$,
\[
  \int f(y,g)\gx (g) dm_G(g) = 0
\]
for a.e $y \in Y$. Then $f$ is orthogonal to the space of $T_{\gr}$
invariant functions. 
\end{pro}

\begin{proof} \cite{FuW96} lemma $9.2$.
\end{proof}

\begin{ntt}
Denote  $U_d$ = $d$ dimensional unitary matrices, $C(U_d)$ the center of
$U_d$ (scalar matrices), and 
$P:U_d \ra \BP U_d = U_d/C(U_d)$ the natural projection. For $U,V \in U_d$
denote by $[U,V]$ the commutator of $U,V$; i.e., $[U,V]=UVU^{-1}V^{-1}$.
\end{ntt}

We need the following lemma:
\begin{lma}\label{abelian}
Let $H$ be a compact abelian connected group, and  
$A:H \ra U_d $ a measurable function.
If $P \circ A$ is a homomorphism, then 
$A(H)$ is a commuting set of matrices.
\end{lma}

\begin{proof}
Let $g,h \in H$. Suppose $[A(h),A(g)]=\gd I$.
If $v$ is an eigenvector of $A(h)$ with eigenvalue $\gc$, then
\[
  A(h)A(g)v=\gd A(g)A(h)v = \gc \gd A(g)v 
\]
thus $A(g)v$ is an eigenvector of $A(h)$ with eigenvalue $\gc \gd$.
This implies that  $A(g)^kv$ is an  eigenvector of $A(h)$ with 
eigenvalue $\gc \gd^k$, thus $\gd$ is a root of unity of order $\le d$.
Denote $C_{d!}$ the group of order $d!$ roots of unity.
Then the commutator set 
\[
  \{ [A(h),A(g)] \}_{h,g \in H} \subset C_{d!}I,
\]
Fix $g$. The function $h \ra [A(h),A(g)]$ is a measurable homomorphism
to  $C_{d!}$ 
\[ \begin{aligned}
  {[A(h_1+h_2),A(g)]} &=[cA(h_1)A(h_2),A(g)]\\
  &=[A(h_1),A(g)][A(h_2),A(g)],
\end{aligned}\]
therefore continuous, and as $H$ is connected it is trivial.
\end{proof}

\begin{thm}\label{sthmCLhom}
Let ${Y}$ be an ergodic m.p.s, and let 
${W}={Y} \times_{\gr} H$ be an ergodic 
extension by a connected abelian group.
Let   $F:Y^0 \times H \times H \ra S^1$  be a measurable function.
Let  $\gs_1(y,h_1)$,$\gs_2(y,h_2):Y^0 \times H \ra S^1$
be measurable functions.
Suppose
\[
 \gs_1(y,h_1)\gs_2(y,h_2)
 = \frac{F(Ty,h_1+\gr(y),h_2+\gr(y))}{F(y,h_1,h_2)}.
\]     
Then for $i=1,2$ there exist measurable functions $g_i,G_i :Y^0 \ra S^1$ 
such that   
\[
 \gs_i(y,h)=g_i(y)\frac{G_i(T_W(y,h))}{G_i(y,h)}
\] 
\end{thm}

\begin{proof}
We construct the following systems: 
for $i=1,2$ let ${X}_i={W} \times_{\gs_i} S^1$, and 
${X}={X}_1 \times_{{Y}} {X}_2$.
Then $\mu_X$
is defined as the conditional product measure relative to the diagonal measure on
$Y^0 \times Y^0$.
The function 
\begin{equation}\label{F_tilde}
 \ti{F}(y,h_1,h_2,\zeta_1,\zeta_2)= 
     F(y,h_1,h_2)\zeta_1^{-1}\zeta_2^{-1}
\end{equation}
is invariant under $T_X$, and therefore by Theorem \ref{max_isometric}
it is measurable with respect to 
$\hat{{Y}}_1 \times \hat{{Y}}_2$ , where $\hat{{Y}}_i$ is the maximal isometric
extension of ${Y}$ in ${X}_i$ for $i=1,2$. 
Isometric extensions are spanned by finite
rank modules (see \ref{isometric_extensions}). Thus 
\[
  \ti{F}(y,h_1,h_2,\zeta_1,\zeta_2)
  = \sum \ip{\ve{\gp^1_j}(y,h_1,\zeta_1)}{\ve{\gp^2_j}(y,h_2,\zeta_2)}
\]  
where 
\[ \begin{aligned}
  T_{X_1}\ve{\gp^1_j}(y,h_1,\zeta_1)
  &=U^1_j(y)\ve{\gp^1_j}(y,h_1,\zeta_1)\\
  T_{X_2}\ve{\gp^2_j}(y,h_2,\zeta_i)
  &=U^2_j(y)\ve{\gp^2_j}(y,h_2,\zeta_2),
  \end{aligned}
\]
and $U^1_j(y),U^2_j(y)$ are  $d_j\times d_j$ unitary matrices. Substituting 
the Fourier expansions:
\[\begin{aligned}
  \ve{\gp^1_j}(y,h_1,\zeta_1)&=\sum \ve{\gp^1_{j,k}}(y,h_1)\zeta^k_1 \\
  \ve{\gp^2_j}(y,h_2,\zeta_2)&=\sum \ve{\gp^2_{j,k}}(y,h_2)\zeta^k_2
 \end{aligned}
\]
in equation(\ref{F_tilde}) we get
that for $k=-1$ there exists $j$ such that $\ve{\gp^1_{j,-1}} \ne 0$.
Apply $T_{X_1}$ to get
\[ 
  \gs_1^{-1}(y,h_1)\ve{\gp^1_{j,-1}}(T_W(y,h_1))
  =U^1_j(y)\ve{\gp^1_{j,-1}}(y,h_1)
\]
For simplicity we drop the indices:
\begin{equation}\label{sMateq}
 \gs^{-1}(y,h)\ve{\gp}(T_W(y,h)) =U(y)\ve{\gp}(y,h)
\end{equation}
For each $y$ consider the distribution of $\ve{\gp}(y,h)$ in the 
fiber over $y$, and look at the vector space spanned by the support 
of this distribution. Call this $V_y$, so that $V_y \subset \BC^d$,
and $V_{Ty}=U(y)V_y$. Since $U(y)$ is unitary, 
$dim V_{Ty}=dim V_y$, thus by ergodicity $dim V_y=\hat{d}$
for a.s. $y$. For each $y$ choose a basis for $\BC^d$ s.t. $V_y$ 
is spanned by the first $\hat{d}$ elements. As the transformation matrix
is a function of $y$, we may assume $d=\hat{d}$.

Denote by $\ve{\ti{\gp}}$ the projection of $\ve{\gp}$ on $\BP V$, and
by $\ti{U}$ the projection of $U$ on $\BP U_d$. Thus:
\[
 \ve{\ti{\gp}}(T_W(y,h))= \ti{U}(y)\ve{\ti{\gp}}(y,h).
\]
Consider the group extension ${W} \times_{\ti{U}} \BP U_d$. 
Then
\[
  \ve{\ti{\gp}}(T_W^n(y,h))= \ti{U}^{(n)}(y)\ve{\ti{\gp}}(y,h)
\]
For fixed $y$, $\{\ve{\ti{\gp}}(y,h)\}_{h \in H}$ spans the space, so 
whenever  $T_W^n(y,h)=(T^ny,h+\gr^{(n)}(y))$ is close to $(y,h)$ 
(by ergodicity this happens for
a generic $y$, and is independent of $h$), $ \ti{U}^{(n)}(y)$ is close to the
identity. This implies that the foregoing group 
extension is not ergodic, and furthermore
- the Mackey group is trivial. Thus for some projective unitary matrix 
function $\ti{M}$:
\begin{equation}\label{M}
 \ti{M}(T_W(y,h))=\ti{U}(y)\ti{M}(y,h)
\end{equation}
Also for any $h'$
\[
 \ti{M}(T_W(y,h+h'))=\ti{U}(y)\ti{M}(y,h+h').
\]
Thus
\[
  \ti{M}^{-1}(T_W(y,h+h'))\ti{M}(T_W(y,h))
 =\ti{M}^{-1}(y,h+h')\ti{M}(y,h)
\]
By ergodicity
\[
  \ti{M}^{-1}(y,h+h')\ti{M}(y,h)=\ti{A}(h'),
\]
for all $h'$, a.e.$(y,h)$.
By Fubini's theorem there exists $h_0$ such that
\begin{equation}\label{phi_0}
  \ti{M}(y,h)=\ti{M}(y,h_0)\ti{A}^{-1}(h-h_0),
\end{equation}
for a.e.$(y,h)$.
The function $\ti{A}(h')$ is a homomorphism of $H$:
\[
 \begin{aligned}
  \ti{A}(h'+h'') &=\ti{M}^{-1}(y,h+h'+h'')\ti{M}(y,h)     \\
                     &=\ti{M}^{-1}(y,h+h'+h'')\ti{M}(y,h+h')
                       \ti{M}^{-1}(y,h+h')\ti{M}(y,h)           \\
                     &= \ti{A}(h'') \ti{A}(h')
 \end{aligned}
\]
Recall $P:U_d \ra \BP U_d$ is  the natural projection. We can find a 
measurable function $A:H \ra U_d$ so that $P\circ A=\ti{A}$. 
\[
A(H) \subset P^{-1}\ti{A}(H).
\]
 Then by lemma \ref{abelian},
$A(H)$ is a commuting set.
Substituting equation  (\ref{phi_0}) in equation (\ref{M}) we get 
\[
 \begin{aligned}
  \ti{M}(Ty,h_0)\ti{A}^{-1}(h+\gr(y)-h_0)
   =& \ti{M}(Ty,h+\gr(y))\\
   =& \ti{U}(y)\ti{M}(y,h)\\
   =& \ti{U}(y)\ti{M}(y,h_0)\ti{A}^{-1}(h-h_0)\\
 \end{aligned}
\]
Thus 
\[
 \ti{U}(y)=\ti{M}(Ty,h_0)\ti{A}^{-1}(\gr(y))\ti{M}^{-1}(y,h_0)
\]
or
\begin{equation}\label{H-A}
  U(y)=M(Ty,h_0)A(-\gr(y))
              M^{-1}(y,h_0)d(y)
\end{equation}
where $d(y)$ is a scalar matrix.
As $A(H)$ is a commuting set, it is simultaneously diagonalizable:
\begin{equation}\label{diagonizable}
  A(h)=N^{-1}D(h)N
\end{equation}
therefore
\[
  U(y)=M(Ty,h_0)N^{-1}D(-\gr(y))N
              M^{-1}(y,h_0)d(y)
\]  
Denote $M'(y)=M(y,h_0)$.
Substitute $U(y)$ in equation (\ref{sMateq}):
\[
 \begin{aligned}
  \gs^{-1}(y,h)NM'^{-1}(Ty) & \ve{\gp}(Ty,h+\gr(y))
  =D(-\gr(y))d(y)NM'^{-1}(y)\ve{\gp}(y,h)
 \end{aligned}
\]
Now each coordinate gives us the desired result.
\end{proof}

\begin{rmr}\label{H_not_connected}
If $H$ in Theorem \ref{sthmCLhom} is not necessarily connected, but
the cocycle $\gr$ is cohomologous to a constant: 
$\gr(y)=c\frac{f(Ty)}{f(y)}$, then we do not need to use lemma \ref{abelian},
and the result holds 
as for some scalar matrix $d(y):Y^0 \ra S^1$
\[
 A(\gr(y))=A(cf(Ty)f^{-1}(y))=A(f(Ty))A(c)A^{-1}(f(y))d(y)
\]
Now diagonalize $A(c)$ : $A(c)=UDU^{-1}$ and substitute in equation
(\ref{H-A}).
\end{rmr}

\begin{thm}\label{thm:hom_general}
Let ${Y}=(Y^0,\caB_Y,\mu_Y,T_Y)$ 
be an ergodic m.p.s. For $i=1,\ldots,k$ 
let ${Y}_i=(Y^0,\caB_Y,\mu_Y,T_Y^i)$. Let 
${W}={Y} \times_{\gr} H$ be an ergodic group extension, where $H$ is a connected abelian group, and let ${W}_i={Y}_i \times_{\gr^{(i)}} H$ 
(notice that $T_{W_i}=T_{W}^i$). 
Let $\nu$ be a joining of the measures on $Y_i$ and,
let $\mu$ be a measure on $\Pi W^0_i$ that is the conditional product 
measure with respect to $\nu$.
Let   $F:\Pi W^0_i \ra S^1$  be a $\mu$ measurable function. 
For $i=1,\ldots,k$, let  $\gs_{i}:W^0 \ra S^1$
be measurable functions, and denote $\pi:W^0 \ra Y^0$ the projection.
Suppose $\mu$ a.e.
\[
 \p_{i=1}^k\gs_i^{(i)}(w_i)
 = \frac{F(T_{W_1}w_1,\ldots,T_{W_k}w_k)}
        {F(w_1,\ldots,w_k)}.
\]     
Then there exist measurable functions $g_i,G_i :Y^0 \ra S^1$ such that   
\begin{equation}\label{eq:W_i}
 \gs_i^{(i)}(w)=g_i(\pi(w))\frac{G_i(T_{W_i}(w))}{G_i(w)}
\end{equation}
\end{thm}
\begin{proof}
The proof is similar to the proof of Theorem \ref{sthmCLhom}.
For $i=1,\ldots,k$ let $X_i=W_i \times_{\gs_i} S^1$. Let
$X$ be the system with $X^0=\p {X}^0_i$, $\mu_X$ the conditional product measure 
with respect to $\nu$, and $T_X=T_{X_1}\times \ldots \times T_{X_k}$.
The function 
\begin{equation}
 \ti{F}(w_1,\gz_1,\ldots,w_k,\gz_k)= 
     F(w_1,\ldots,w_k)\zeta_1^{-1}\ldots \zeta_k^{-1}
\end{equation}
is invariant under $T_X$. Proceeding  as in Theorem \ref{sthmCLhom}.
we find that equation (\ref{eq:W_i}) holds 
on the ergodic components of $T_W^i$. As $T_W$ is ergodic, $T_W^i$ has finitely
many ergodic components. Let $Y_l$ be an ergodic component of  $T_Y^i$. The ergodic components of $T_W^i$ which project onto $Y_l$ are determined by the Mackey 
group $M_l$ which is a closed subgroup of $H$. As  $T_W^i$ has finitely many ergodic 
components $M_l$ is of finite index in $H$, but 
$H$ is connected, therefore has no closed subgroups of finite index.  
Therefore the  ergodic components of $T_W^i$ are  of the form $Y_l \times H$.
\end{proof}

\begin{cor}\label{hom} Let ${Y}$ be an ergodic m.p.s, 
${X}={Y} \times_{\gr} H$ an ergodic abelian 
extension where either  $H$ is connected or  
the cocycle $\gr$ is cohomologous to a 
constant.
Suppose there exists a measurable family of measurable functions 
$\{ f_u\}_{u \in H}$,
$f_u:Y^0 \times H \ra S^1$ such that 
\[
\frac{\gs(y,h+u)}{\gs(y,h)} 
= \frac{f_u(T_X(y,h))}{f_u(y,h)},
\]
then there exist measurable functions $g:Y^0 \ra S^1$ and 
$F:Y^0 \times H \ra S^1$ such that 
\[
 \gs(y,h)=g(y)\frac{G(T_X(y,h))}{G(y,h)}.
\]
\end{cor}

\begin{proof}
Make the coordinate change: $h_1=h; \ h_2=h+u$.
Then 
\[
  f_u(y,h)=f(y,u,h)=f'(y,h+u,h)=f'(y,h_1,h_2)
\]
and
\[
  f_u(Ty,h+\gr(y))=f'(Ty,h_1+\gr(y),h_2+\gr(y)).
\]
Now apply Theorem \ref{sthmCLhom}.
\end{proof}

\begin{lma}\label{ef_times_char} Let ${Y}={Z} \times_{\gr} H$ be an ergodic 
abelian extension of ${Z}$, and 
$F: Z^0 \times H \ra S^1$, $g: Z^0 \ra S^1$ measurable functions such that
\[
  g(z)=\frac{T_{Y}F(z,h)}{F(z,h)}.
\]
Then there exists $\gx \in \hat{H}$, and $k:Z^0 \ra S^1$ such that
\[
  F(z,h)=k(z)\gx(h).
\]
\end{lma}

\begin{proof}
Take the Fourier expansion of $F$:
\[
  F(z,h)= \sum k_i(z)\gx_i(h).
\]
Then for all $i$
\[
 k_i(T_Zz)\gx_i(h)\gx_i(\gr(z))=g(z)k_i(z)\gx_i(h).
\]
Ergodicity of $T_Z$ implies $| k_i(z) |$ is constant a.e.
The fact that $|F|=1$ implies that there exist an $i$ for which
$|k_i(z)| \ne 0$. If there are two such indices $i,j$, then
\[
 \frac{\gx_i}{\gx_j}(\gr(z))
\] 
is a coboundary. As $T_{Y}$ is ergodic $\gx_i/\gx_j=1$
(otherwise the Mackey group of the extension $Z\times_{\gr}H$ is not $H$).
\end{proof}

\begin{ntt} Let $(X^0_1,\caB_1)$, $(X^0_2,\caB_2)$ be measure spaces.
Denote
\[
 B(X^0_1,X^0_2)=\{f:X^0_1 \ra X^0_2, f \ \te{measurable} \}.
\]
\end{ntt}

\begin{lma}\label{measurable}
Let ${Y}={Z} \times_{\gr} H$ be an ergodic abelian extension of ${Z}$,

$(X,\mu)$ a measure space, and let
$x \ra f_{x}(y)$ be a Borel measurable function
from $X$ to $B(Y^0,S^1)$. 
Suppose for all $x \in X$ there are functions
$g_{x}(z),F_{x}(y)\in B(Y^0,S^1)$ such that
\begin{equation}\label{measurability}
  f_{x}(y)=g_{x}(z)\frac{T_YF_{x}(y)}{F_{x}(y)}.
\end{equation}
Then there is a $\mu$ measurable choice of $g_{x}(z),F_{x}(y)$.
\end{lma}

\begin{proof}
Endowed with the $L^2$ topology, $B(Y^0,S^1)$ is a polish group. 
Let $B(Z^0,S^1)$ be the closed subgroup of $B(Y^0,S^1)$ of functions 
that depend only on the $z$ coordinate, and let $f \ra \bar{f}$ be the 
natural projection onto $\bar{B}=B(Y^0,S^1)/ B(Z^0,S^1)$, with the induced
topology. By a theorem of Dixmier (\cite{BK96} Theorem $1.2.4$)
there is a measurable section $\bar{B} \ra B$.  
Equation (\ref{measurability}) implies
\[
  \bar{f}_{x}(y)=\frac{T_Y\bar{F}_{x}(y)}{\bar{F}_{x}(y)}.
\] 
Define $\gf:\bar{B} \ra \bar{B}$
\[
 \gf(\bar{f})=\frac{T_Y\bar{f}}{\bar{f}}.
\]
If $\gf(\bar{f})=\gf(\bar{g})$, then for some function $h(z)$
\[
  \frac{T_Y\frac{f}{g}(y)}{\frac{f}{g}(y)}= h(z).
\]
By \ref{ef_times_char} this implies that up to multiplication by 
a function of $z$,
$\frac{f}{g}$ belongs to a countable set, thus $\gf$ is countable to one.
By Lusin \cite{Lu30} 
$\gf(\bar{B})$ is a measurable set and there is a measurable 
function $\gp: \gf(\bar{B}) \ra \bar{B}$ s.t. 
\[
 \gf \circ \gp =Id|_{\gf(\bar{B})}
\]
Now if 
\[
 \gp(\bar{f}_{x})=\bar{F}_{x},
\]
then
\[
 \bar{f}_{x}=\gf \circ \gp(\bar{f}_{x})
              =\frac{T_Y\bar{F}_{x}}{\bar{F}_{x}}.
\]
The composition
\[
 x \ra f_{x} \ra \bar{f}_{x} \ra \bar{F}_{x} \ra F_{x}
\]
gives a measurable choice of $F_{x}$, and $g_{x}$ is measurable as
a quotient of measurable functions.
\end{proof}

\begin{rmr}\label{g_x_constant}
If $g_x(z) \in B(Y^0,*)$ ($g_x$ is constant) then
the same proof works to give a measurable choice of $g_x,F_x$. 
\end{rmr}

\begin{ntt}
We write  $f \sim  g$ if $f/g=$const.
\end{ntt}

\begin{lma}\label{multiplicity_of_f_h}
Let ${X}={Y} \times_{\gr} H$ be an ergodic abelian extension of ${Y}$.
Let $\gs :Y^0 \times H \ra S^1$ be such that for all $u \in H$ there 
exists a measurable function $f_u : Y^0 \times H \ra S^1$ and a constant 
$\gl_u$ such that
\begin{equation}\label{quazi}
  \frac{\gs(y,h+u)}{\gs(y,h)}=\gl_u \frac{f_u(T_X(y,h))}{f_u(y,h)},
\end{equation}
Then there exists a measurable family of measurable functions 
$\{ f_{u} \}_{u \in H}$, a measurable family of constants
$\{ \gl_{u} \}_{u \in H}$
satisfying the above equation, and a  neighborhood  of zero $U$ 
in $H$ such that 
\[ \begin{aligned}
 f_{u_1+u_2}(y,h) &\sim f_{u_2}(y,h+u_1)f_{u_1}(y,h) \\
 \gl_{u_1+u_2}&=\gl_{u_1}\gl_{u_2} \end{aligned}
\]
whenever $u_1,u_2, u_1+u_2 \in U$.
\end{lma}

\begin{proof}
By remark \ref{g_x_constant} we may assume that the families 
$\{ f_{u} \}_{u \in H}$, and $\{ \gl_{u} \}_{u \in H}$ depend 
measurably on $u$.
Using equation (\ref{quazi}) we get
\[
 \begin{aligned}
  \frac{\gs(y,h+u_1+u_2)}{\gs(y,h)}
  &= \gl_{u_1+u_2}
     \frac{T_Xf_{u_1+u_2}(y,h)}
          {f_{u_1+u_2}(y,h)} \\
  &= \gl_{u_1}\gl_{u_2}
    \frac{T_Xf_{u_1}(y,h+u_2)}{f_{u_1}(y,h+u_2)}
    \frac{T_Xf_{u_2}(y,h)}{f_{u_2}(y,h)} 
\end{aligned} 
\]
this implies that 
\[
 \frac{f_{u_1+u_2}(y,h)}{f_{u_1}(y,h+u_2)f_{u_2}(y,h)}
\]
is an eigenfunction of $T_X$ and that 
\[
 \frac{\gl_{u_1}\gl_{u_2}}{\gl_{u_1+u_2}}
\] 
is an eigenvalue. Let ${Z}$ be the Kronecker factor of ${X}$, 
$\pi :X^0 \ra Z^0$ the projection map, let $N$ parametrize $\hat{Z}$,
and let $\gp_{N(u_1,u_2)}(z)$ be a character of $Z$ s.t.:
\begin{equation}\label{eigenfn}
  \frac{f_{u_1+u_2}(y,h)}{f_{u_2}(y,h+u_1)f_{u_1}(y,h)}
  \sim \gp_{N(u_1,u_2)}\circ \pi(y,h)
\end{equation}
and 
\begin{equation}\label{eigenva}
  \frac{\gl_{u_1}\gl_{u_2}}{\gl_{u_1+u_2}}
  =\gp_{N(u_1,u_2)}(\ga).
\end{equation}
Any two characters taking the same value on $\ga$ are the same, therefore
$\gp_{N(u_1,u_2)}$ is symmetric, i.e
\[\gp_{N(u_1,u_2)}=\gp_{N(u_2,u_1)}\]
We now show that $\gp_{N(u_1,u_2)}$ satisfies a $2$-cocycle equation:
\[
 \begin{aligned}
  &\gp_{N(u_1+u_2,u_3)}\circ \pi (y,h)   \sim
   \frac{f_{u_1+u_2+u_3}(y,h)}
        {f_{u_3}(y,h+u_1+u_2)f_{u_1+u_2}(y,h)} \\
  &\gp_{ N(u_1,u_2+u_3)}\circ \pi(y,h) \sim
    \frac{f_{u_1+u_2+u_3}(y,h)}
         {f_{u_2+u_3}(y,h+u_1)f_{u_1}(y,h)}
  \end{aligned}
\]
Thus
\[
 \begin{aligned}
   (\gp_{N(u_1+u_2,u_3)}\circ \pi (y,h))& 
   f_{u_3}(y,h+u_1+u_2)f_{u_1+u_2}(y,h)  \\
  \sim & (\gp_{N(u_1,u_2+u_3)}\circ \pi (y,h))
    f_{u_2+u_3}(y,h+u_1)f_{u_1}(y,h)
 \end{aligned}
\]
Dividing both sides by 
\[
f_{u_1}(y,h)f_{u_3}(y,h+u_1+u_2)f_{u_2}(y,h+u_1)
\]
we get
\[
 \begin{aligned}
  (\gp_{N(u_1+u_2,u_3)}\circ \pi (y,h))& 
  \frac{f_{u_1+u_2}(y,h)}
       {f_{u_1}(y,h)f_{u_2}(y,h+u_1)} \\
  \sim (\gp_{N(u_1,u_2+u_3)} & \circ \pi (y,h))
  \frac{f_{u_2+u_3}(y,h+u_1)}
       {f_{u_2}(y,h+u_1)f_{u_3}(y,h+u_1+u_2)}.
 \end{aligned}
\]
Combining the above equation with equation (\ref{eigenfn}),
\begin{equation}\label{2cocycleN}
 \begin{aligned}
  \gp_{N(u_1+u_2,u_3)}\gp_{N(u_1,u_2)} 
  =
  \gp_{N(u_1,u_2+u_3)} \gp_{N(u_2,u_3)}.
 \end{aligned}  
\end{equation}

As $u \ra f_{u}$  is a measurable 
function, $f_{u_2}(y), f_{u_2+u_1}(y)$ are close in measure for small
$u_1$, most $u_2$, and the same goes for  
$f_{u_2}(y,h) , f_{u_2} (y,h+u_1)$. Therefore the expression
in equation (\ref{eigenfn}) is close (in measure) to 
$\bar{f_{u_1}}(y,h)$.
But $N_1 \ne N_2$ implies 
\[
  \|\gp_{N_1}-\gp_{N_2}\|_2 =\sqrt{2},
\]
thus by  equation (\ref{eigenfn}),
$\gp_{N(u_1,u_2)} = \gp_{\ti{N}(u_1)}$
for $u_1 \in U'$ a neighborhood of zero in $H$, 
$u_2 \in A$ a set of positive measure. 
The set  $A-A$ contains a neighborhood
of zero $U''$. Let $U=U' \cap U''$.
Take any $u_1,u_2, u_1+u_2 \in U$, and find an element 
$u_3 \in A$ such that 
$u_3+u_2 \in A$, then by (\ref{2cocycleN})
\[
 \gp_{N(u_1,u_2)}= \gp_{\ti{N}(u_1)}
    \gp^{-1}_{\ti{N}(u_1+u_2)}\gp_{\ti{N}(u_2)}
\]
For $u  \in U$,
denote 
\[ 
  \ti{f}_{u}(y,h) = (\gp_{\ti{N}(u)}\circ \pi (y,h)) f_{u}(y,h),
\]
and 
\[
  \ti{\gl}_{u}=\gl_{u}\gp^{-1}_{\ti{N}(u)}(\ga).
\]
By equations (\ref{eigenfn}), if 
$u_1, u_2,u_1+u_2 \in U$, then : 
\begin{equation}\label{eigenfn1}
  \ti{f}_{u_1+u_2}(y,h) \sim \ti{f}_{u_2}(y,h+u_1)
  \ti{f}_{u_1}(y,h).
\end{equation}
By equation(\ref{eigenva}),
if $u_1, u_2,u_1+u_2 \in U$ then
\begin{equation}\label{eigenva1}
 \ti{\gl}_{u_1+u_2}=\ti{\gl}_{u_1}\ti{\gl}_{u_2}.
\end{equation}
\end{proof}

\begin{lma}\label{reduction_finite_connected}
Let $H$ be a torus (possibly infinite dimensional) 
and let ${X}={Y} \times_{\gr} H$ be an ergodic abelian extension of ${Y}$.
Suppose
\begin{equation}
  \frac{\gs(y,h+u)}{\gs(y,h)}=\gl_u \frac{f_u(T_X(y,h))}{f_u(y,h)}
\end{equation}
for each $u \in H$, and $\gl_u$ and $f_u$ depend measurably on $u$.
Then there is a  subgroup $J < H$ such that  $H/J=\BT^n$ such that 
if $\pi: H \ra H/J$ is the natural projection then
there exists a function $\ti{\gs}:Y^0 \times (H/J) \ra S^1$ such that 
\[
  \gs(y,h)=\ti{\gs}(y,\pi(h))\frac{F(T_{X}(y,h))}{F(y,h)}.
\]
\end{lma}

\begin{proof}
By lemma \ref{multiplicity_of_f_h} the functions $f_u$ can be chosen 
such that $\gl_u$ is multiplicative in  a zero neighborhood $U$
in $H$. The neighborhood  $U$ contains $J_1$ - 
a closed connected subgroup of H, such that $H/J_1=\BT^l$,
thus $\gl_h$ is a character of $J_1$. 
Thinking of $H$ (measurably) as $H/J_1 \times J_1$ with coordinates $(h_0,j)$
the above equation becomes
\begin{equation}
  \frac{\gs(y,h_0,j+u)}{\gs(y,h)}=
  \gl_u \frac{f_u(Ty,h+\gr(y))}{f_u(y,h)}.
\end{equation}
where $u \in J_1$. 
This is the same as
\begin{equation}
  \frac{\gl^{-1}_{j+u}\gs(y,h_0,j+u)}{\gl^{-1}_{j}\gs(y,h)}=
   \frac{f_u(Ty,h+\gr(y))}{f_u(y,h)}.
\end{equation} 
Applying corollary \ref{hom} replacing $Y$ with  $(Y \times H/J_1)$
and $H$ with $J_1$ we get
\[
  \gl^{-1}_{j}\gs(y,h_0,j)=\ti{\gs}(y,h_0)\frac{T_{\gr}F(y,h)}{F(y,h)}.
\]  
or
\[
  \gs(y,h_0,j)=\gl_{j}\ti{\gs}(y,h_0)\frac{T_{\gr}F(y,h)}{F(y,h)}.
\] 
Now for $j$ in the kernel of $\gl$ we have $\gl_{j}=1$. The image
of $\gl$ is $S^1$ thus if $ker \gl$ is $J$ then $H/J=\BT^{l+1}$.

\end{proof}

\begin{rmr}\label{rmr:co-tori} 
The group $H$ is a compact connected abelian (metrizable) group, 
and therefore has countably many closed subgroups $J$
such that $H/J$ is a finite dimensional torus. 
\end{rmr}

\begin{rmr}\label{reduction_no_torus}
If $H$ is any connected compact abelian group (not necessarily a torus)
then $J_1$ in the foregoing proof is not necessarily connected. 
By the same proof we will get that $\gs$ is cohomologous to a 
cocycle lifted from a product of a finite torus and a totally disconnected
compact abelian group.
\end{rmr}

\begin{lma}\label{reduction_finite_disconnected}
Let ${X}={Y} \times_{\gr} H$ be an abelian extension of ${Y}$ with $\gr(y)$ 
cohomologous to a constant function (now $H$ is any compact abelian group, 
not necessarily connected). 
Let $\gs$ be as in lemma \ref{reduction_finite_connected}.
Then there is a  subgroup $J < H$, and a finite group $C_k$,
such that  $H/J=\BT^n \times C_k$ 
and letting $\pi: H \ra H/J$ denote the natural projection then
there exists a function $\ti{\gs}:Y^0 \times H/J \ra S^1$ such that 
\[
  \gs(y,h)=\ti{\gs}(y,\pi(h))\frac{T_{X}F(y,h)}{F(y,h)}.
\]
\end{lma}

\begin{proof}
By lemma \ref{multiplicity_of_f_h} the functions $f_u$ can be chosen such 
that $\gl_u$ is multiplicative in  a zero neighborhood $U$
in $H$. The neighborhood  $U$ contains $J_1$ - a closed subgroup of $H$,
such that $H/J_1=\BT^l \times C_j$ where $C_j$ is a finite group,
thus $\gl_u$ is a character of $J_1$. Now proceed as in lemma 
\ref{reduction_finite_connected} (the image of $\gl$ is either $S^1$ or a finite group).  
\end{proof}

\section{Lie Groups and Nilsystems}

\begin{dsc}\label{dsc:ergodic_com}
Let $N$ be a $k$-step, simply connected nilpotent Lie group,  
$\gC$ a discrete subgroup s.t.
$N/\gC$ is compact. Let $\caB$ be the (completed) Borel algebra, 
$m$ the Haar measure  on $N/\gC$, and let $a \in N$.
The system ${X}=(N/\gC,\caB,m,T)$ where $Tg\gC=ag\gC$ is called 
a {\em nilflow}. We will sometime denote this system $(N/\gC,a)$. 
Let $N_1=N$, and for $i>1$: $N_i=[N_{i-1},N]$ 
($N$ is a $k$-step nilpotent group if $N_{k+1}=\{1\}$), and for $i \ge 1$
let  $\gC_i=\gC \cap N_i$. The groups $N_j$ for $j>1$ are connected 
(see \cite{L02}).
Th group $N_k$ is abelian and connected and therefore isomorphic to 
$\BR^m$ for some positive integer $m$. 
Then $\gC_i$ is a discrete subgroup
of $N_i$, and $N_i/\gC_i$ is compact.
Let $m_i$ be the Haar measure on $N_i/\gC_i$.
Let
\[
M=\{(y_1,y_1^2y_2,\ldots, \p_{j=1}^{k} y_j^{\binom{l}{j}}):y_1 \in N_1, 
     y_2 \in N_2,\ldots,y_{k} \in N_{k}\} \subset N^{l},
\]
where 
\[
\p_{j=1}^{k} y_j^{\binom{l}{j}}=y_1^k y_2^{\binom{l}{2}} \ldots
           y_k^{\binom{l}{k}}.
\]  
The elements of $M$ are called 
Hall-Petresco (HP) sequences, and form a group (see \cite{La54}, \cite{L98}).
The first $k+1$ elements in the sequence determine the rest. 
Computation 
shows that if $n \in N$ and $(n_1, \dots,n_l) \in M$ then 
$([n,n_1],\ldots,[n,n_l]) \in M$.
Let $\gL=M \cap \gC^{l}$. The nilmanifold $Y=M/\gL$ is 
embedded in $(N/\gC)^{l}$ and let $\nu$ be the Haar measure on $Y$.
Then for almost all $g \in N$ for all $(n_1,\ldots n_l) \in M$ 
\begin{equation}\label{eq:avr_on_nil}
   \Lim{N}  \Avr{N}{n} \p_{j=1}^l f_j(a^{jn}n_jg\gC)
         =\int_Y \p_{j=1}^l f_j(gz_j)          
             d\nu(z_1,\ldots,z_{l}).
\end{equation}
For more details see \cite{Z02a}. A similar result holds for the action
of \\$(a^{m_1},\ldots ,a^{m_l})$ for any $m_i \in \BZ$.
\end{dsc}

\begin{rmr}\label{gcd}
If $N/\gC$ is connected then $a^k$ is ergodic for any $k \ne 0$
therefore equation \ref{eq:avr_on_nil}
remains the same if we replace $a$ by $a^k$.
\end{rmr}

\begin{dsc}
Let $(N/\gC,a)$ be a $k$-step nilflow. Define
\[ \begin{split}
   \tau_{l}(a):=& a \times \ldots \times a^{l}\\
   \triangle_{l}(a):=& a \times \ldots \times a.
   \end{split}
\] 
Let $\triangle_l(m)$ be the diagonal measure on $(N/\gC)^l$.\\ 
Define a measure on $(N/\gC)^l$
\begin{equation}\label{def_mu_k}
 \bar{\triangle}_{l}(m):= 
 \lim_{N \ra \infty} \Avr{N}{n} \tau_{l}(a)^n \triangle_{l}(m).
\end{equation}
By the above discussion the ergodic components of
$\bar{\triangle}_{l}(m)$ are parametrized by $N/\gC$ and are
of the form $M(g\gC,\ldots,g\gC)$.
\end{dsc}

\begin{dsc}\label{abelian_ext_nil} 
The system $(N/\gC,a)$ may be represented as an Abelian extension
of a $k-1$-step nilflow 
$(N/\gC,a)=(N/N_k\gC \times_{\gr} N_k/\gC_k)$. Inductively $(N/\gC,a)$
may be represented as a tower of Abelian extensions, starting out with a point.
(the first block in the tower would be the Kronecker factor $N/N_{2}\gC$).
Consider the system $Y=((N/N_k\gC)^l, \bar{\triangle}_{l})$. Then 
$(Y \times_{(\gr^{(1)},\ldots,\gr^{(l)})} (N_k/\gC_k)^l,\tau_l(a))$ is an 
abelian group extension of $Y$. 
The Mackey group, associated with the ergodic components of the group 
extension that are mapped onto the ergodic components of $\bar{\triangle}_{l}$,
is 
\[
M_{k,l}=
\{(g_1,g_1^2g_2,\ldots, \p_{j=1}^{k} g_j^{\binom{l}{j}}):g_1,\ldots,g_{k} \in N_{k}\}\left/ \gC_k^l \right.
\]
In additive notation, denote $H=N_k/\gC_k$, then
\[
M_{k,l}=
\{(h_1,2h_1+h_2,\ldots, \sum_{j=1}^{k} \binom{l}{j}h_j) : h_1,\ldots,h_{j} \in 
H \}.
\]

\begin{lma}\label{Y_j_of_nil}
Let $X=(N/\gC,a)$ be a $k$-step nilflow. Then ${Y}_{r}(X)=(N/N_r\gC,a)$ for 
$r\ge 2$.
\end{lma}

\begin{proof}
We prove this by induction on the nilpotency level  of $N$. Let $N$ be a 
$1$-step nilpotent group (thus $(N/\gC,a)$ is a Kronecker system).
Let $\gp$ be an eigenfunction of $(N/\gC,a)$, then
\[
\Avr{N}{n} T^n (\gp)^2 T^{2n} (\gp)^{-1} =\gp.
\] 
Therefore 
\[
L^2(Y_1) = \te{span} \{ \te{eigenfunctions}\}=L^2(N/\gC).
\]
Now assume the statement for $(k-1)$-step nilflows. Let 
$X=(N/\gC,a)$ be a $k$-step nilflow. By lemma \ref{natural}, for $r<k+1$,
$Y_r(N/N_r\gC,a)$ is a factor of $Y_r(X)$. By the induction hypothesis
$Y_r(N/N_r\gC,a)=(N/N_r\gC,a)$. But the integral in equation 
(\ref{eq:avr_on_nil}) for $l=r$, is a function on $N/\gC$
that is invariant under translation by elements of $N_r$. By lemma
\ref{lma:limit_char}, $L^2({Y}_{r}(X))$ is spanned by these integrals. 
Let  $r=k+1$, and let $f(g\gC) \in L^{\infty}(N/\gC)$. 
We want to show that the function $f$ is in the span of the 
integrals one obtains  in (\ref{eq:avr_on_nil}) (the value of the limit
is the same for any choice of $(n_1,\ldots,n_{k+1}) \in M$).
The element
$g\gC$ is determined by the first $k+1$ elements 
$m_1g\gC,\ldots,m_{k+1}g\gC$ of a HP geometric progression. 
Therefore there exists a function
$F \in L^{\infty}{\bt_{k+1}(m)}$ such that 
\[
F(m_1g\gC,\ldots,m_{k+1}g\gC)=f(g\gC).
\] 
Now as 
\[
F(m_1g\gC,\ldots,m_{k+1}g\gC)=F(am_1g\gC,\ldots,a^{k+1}m_{k+1}g\gC)
\]
we get
\[
\lim_{N \ra \infty} \Avr{N}{n} \tau_{k+1}(a)^n 
      F(m_1g\gC,\ldots,m_{k+1}g\gC)=f(g\gC).
\]
Now approximate $F$ (in $L^2(\bt_{k+1}(m))$) by sums of functions of the type 
$f_1\otimes \ldots \otimes f_{k+1}$.
\end{proof}

\begin{cor}\label{cor:ucf_pronil}
If $Y$ is a $k$-step pro-nilflow then $Y_{k+1}(Y)=Y$. 
\end{cor}
\end{dsc}

\begin{lma}\label{multiplicity_f_s_nil} Let $(N/\gC,a)$ be a $k$-step nilflow.
Let $f \in B(N/\gC,S^1)$.
Let $\{\gl_c\}$ be a family of constants, 
$\{f_c\}_{c \in N_{k}}$ be a family of functions in  $B(N/\gC \ra S^1)$ 
such that 
\begin{equation}\label{eq:center_N_nil}
\frac{f(cy)}{f(y)}=\gl_{c} 
\frac{f_{c}(ay)}{f_{c}(y)}
\end{equation}
for all $c \in N_{k}$. Then   
we can choose $f_c,\gl_c$ 
such that 
\[
  f_{c_1}(c_2y) f_{c_2}(y) \sim f_{c_1c_2}(y).
\]
\end{lma}

\begin{proof}
By lemma \ref{multiplicity_of_f_h} this holds in a neighborhood
of zero $U \subset N_{k}$.
Notice that multiplying $f_{c}$ by a constant does not
affect equation (\ref{eq:center_N_nil}). 
If $c \in N_{k}$, $c= c_1\ldots c_s$, and 
$c_1, \ldots,c_s \in U$, define
\[
  f_{c}(y)=f_{c_1}(y)f_{c_2}(c_1y)\ldots
                   f_{c_s}(c_1 \ldots c_{s-1}y).
\]
We claim this is well defined (up to a constant multiple) on 
$N_{k}$:
given two sequences  $c_1,\ldots c_s$ and $c'_1,\ldots c'_t$ with 
equal product, we can break up the `steps' $c_r$ into an equal number of 
small steps and we can interpolate a sequence of such paths where two 
consecutive paths differ only within a small cube which can be translated to 
be inside $U$. Since the resulting $\gl$'s and $f$'s will be the same for 
consecutive paths, they will be the same for the initial and the final ones.
\end{proof}

\begin{lma}\label{matrix_nil} Assume $N/\gC$ has no non trivial finite factors.
The group $N/N_2\cong \BZ \times \BR^n$ for some integer $n$ . 
The action of $a$ on it is given by rotation by some element $(1,\ga)$. 
Under the conditions of lemma 
\ref{multiplicity_f_s_nil} we can choose $f_c,\gl_c$ 
\[
\frac{f(cy)}{f(y)} =e^{2 \pi i \ip{L\ga}{c}}\frac{f_{c}(ay)}{f_{c}(y)},
\]
for an integer $n \times m$ integer matrix $L$. 
\end{lma}

\begin{proof}
Now $\gl_{c}$ is a continuous 
multiplicative function on $N_{k}=\BR^m$.
We now use additive notation for
$N_{k}$. 
In this notation $\gl_{c}$ is  of the
form $e^{2 \pi i \ip{r}{c}}$. Let $e_i$ denote the standard basis for $\BR^m$.
Each $f_{e_i}$ is an eigenfunction (as the left side of equation
(\ref{eq:center_N_nil}) is $1$).
Thus there exists  $\vec{n}_i \in \BZ^n$ such that 
\[
f_{e_i}(y)=Ce^{2 \pi i \ip{\vec{n}_i}{\bar{y}}},
\]
 where $\bar{y} \in N/N_2$,
with eigenvalue $e^{2 \pi i \ip{\vec{n}_i}{\ga}}$. 
Finally for each $i$ there is $k_i \in \BZ$
such that
\[
\ip{r}{e_i}=\ip{\vec{n}_i}{\ga}+k_i.
\]
Now take $L$ the matrix with the $i$'th row being $(k_i,\vec{n}_i)$
(the action of $a$ on $\BZ\times \BR^n$ given by $(1,\ga)$). 
\end{proof}

\begin{lma}\label{commuting_center_nil}
For any $c_1,c_2 \in N_{k}$, 
$f_{c_1}$ and $f_{c_2}$ 
satisfying equation (\ref{eq:center_N_nil}) we have
\[
 \frac{f_{c_1}(c_2y)}{f_{c_1}(y)}=\frac{f_{c_2}(c_1y)}{f_{c_2}(y)}
\] 
\end{lma}

\begin{proof}
The function 
\[
 \frac{f_{c_1}(c_2y)}{f_{c_1}(y)} \left/\frac{f_{c_2}(c_1y)}{f_{c_2}(y)}\right.
\]
is invariant under rotation by $a$, therefore 
\[
 \frac{f_{c_1}(c_2y)}{f_{c_1}(y)}=C(c_1,c_2)\frac{f_{c_2}(c_1y)}{f_{c_2}(y)}
\]
Using lemma \ref{multiplicity_f_s_nil}
\[ \begin{split}
 C(c_1,c_2c) & \frac{f_{c_2}(c_1cy)f_{c}(c_1y)}{f_{c_2}(cy)f_{c}(y)}
 =C(c_1,c_2c)\frac{f_{c_2c}(c_1y)}{f_{c_2c}(y)}
 =\frac{f_{c_1}(c_2cy)}{f_{c_1}(y)}\\
 =&\frac{f_{c_1}(c_2cy)}{f_{c_1}(cy)}\frac{f_{c_1}(cy)}{f_{c_1}(y)}
 =C(c_1,c_2)C(c_1,c)\frac{f_{c_2}(c_1cy)}{f_{c_2}(y)}
  \frac{f_{c}(c_1y)}{f_{c}(y)}
\end{split}
\]
Therefore $C(c_1,c_2)$ is multiplicative in $c_1,c_2$.
If $c \in   N_{k} \cap \gC$ 
then $f_{c}$ is an eigenfunction. Thus for
$c \in   N_{k} \cap\gC$,
$C(c_1,c)=C(c,c_2)$=1. This implies that for any
$c_1 \in N_{k}$, 
$C(c_1,*)$ is a character of 
$N_{k}/(N_k \cap \gC)$ 
which is a compact
connected abelian group. As there are countably many of those 
$C(c_1,c_2) \equiv 1$.   
\end{proof}

\section{the van der Corput lemma} 
One of the main tools in studying characteristic factors is  
the van der Corput lemma. 
The formulation below is due to 
Bergelson \cite{Be87}:
\begin{lma}[van der Corput]\label{VDC}
Let $\{ u_n \}$ be a bounded sequence of vectors in a Hilbert
space $\caH$. Assume that for each $m$ the limit
\[
  \gc_m := \Lim{N} \Avr{N}{n} \ip{u_n}{u_{n+m}}
\]
exists, and
\begin{equation}\label{vdc}
  \Lim{M} \Avr{M}{m} \gc_m=0.
\end{equation}
Then
\[
  \Avr{N}{n} u_n \overset{\caH}{\lra} 0.
\]
\end{lma} 

\begin{proof}
Let $M$ be large enough so that the expression in (\ref{vdc})
 is small.
Let $N$ be large enough with respect to $M$ so that the 
two expressions
\[
   \frac{1}{NM}\Sum{n}{N}\Sum{m}{M} u_{n+m}, \quad \Avr{N}{n} u_n
\] 
are close. 
We have:
\begin{equation*}
 \begin{split}
   \| \frac{1}{NM}\Sum{n}{N}\Sum{m}{M} u_{n+m}
     \|^2 & \le \Avr{N}{n} \| \Avr{M}{m} u_{n+m} \|^2      \\
   & = \frac{1}{NM^2}\Sum{n}{N}\Sum{m_1,m_2}{M} 
       \ip{u_{n+m_1}}{u_{n+m_2}}                          \\
   & \overset{N \ra \infty}{\lra} \frac{1}{M^2}\Sum{m_1,m_2}{M}
     \gc_{m_2-m_1}
 \end{split}
\end{equation*}
which is small.
\end{proof}

\section{proof of Theorem \ref{main}}
Let $X$ be an ergodic m.p.s, and let $Y_j(X)$ be the $j$-u.c.f of $X$, 
and let $\pi_j:X \ra Y_j(X)$ be the factor map. When the context 
is clear we will write  $T$ for $T_X$, and $Y_j$ for $Y_j(X)$.
Let $\vec{a}=(a_1,\ldots,a_l) \in \BZ^l$. We will always assume that
$a_i$ are distinct. Denote 
\[ \begin{split}
   \tau_{\vec{a}}(T):=& T^{a_1} \times \ldots \times T^{a_l}\\
   \triangle_{l}(T):=& T \times \ldots \times T.
   \end{split}
\] 
When the context is clear we will use $\tau_{\vec{a}}$ for 
$\tau_{\vec{a}}(T)$, and $T$ or $T_l$ for  $\triangle_l(T)$.

Let $\triangle_l(\mu_{X})$ be the diagonal measure on $(X^0)^l$.
We will prove theorem \ref{main} inductively, along with a sequence of 
statements (theorem \ref{main} is item (\ref{item:main})). 
\begin{thm} \label{thm:main} 
\begin{enumerate}
\item \label{def_mu_l}Let $\vec{a}=(a_1,\ldots,a_{j+1}) \in \BZ^{j+1}$. 
      The limit
      \begin{equation*}
        \bar{\triangle}_{\vec{a}}(\mu_X):= 
        \lim_{N \ra \infty} \Avr{N}{n} \tau_{\vec{a}}^n \triangle_{j+1}(\mu_X).
      \end{equation*}
      exists. Furthermore $\bar{\triangle}_{\vec{a}}(\mu_{X})$ is the 
      conditional product measure relative to 
      $\bar{\triangle}_{\vec{a}}(\mu_{Y_{j}})$.
\item \label{item:iso} $Y_{j+1}(X)$ is an isometric extension of $Y_j(X)$.
\item \label{def:type_j} Let $X$ be an ergodic m.p.s. 
      Let $l\in \BN$. Let $\mu$ be a measure on $(X^0)^l$,  
      let $\vec{a}=(a_1,\ldots,a_l) \in \BZ^l$, and for 
      $i=1,\ldots,l$, let $f_i\in B(X^0,S^1)$. Recall that 
      $f^{(m)}(x)=f(T^{m-1}x)\ldots f(Tx)f(x$.) 
      We say that $(f_1,f_2, \ldots,f_{l})$ is of {\em type
      $\vec{a}$} w.r.t $\mu$  
      if there exists a $\mu$-measurable function $F$
      taking values in $S^1$, such that
      \[ 
        \p_{i=1}^{l} f_i^{(a_i)}(x_i)
        =\frac{\tau_{\vec{a}} F(x_1,\ldots,x_l)}{F(x_1,\ldots,x_l)}.
      \] 
      Let $H$ be a compact abelian group, $Y$ a $(j-1)$-step pro-nilflow. 
      We say that
      $\gr:Y \ra H$ is of {\em type $j$} if for any character 
      $\gx \in \hat{H}$, there exists a character 
      $\ti{\gx}=(\gx_1,\ldots,\gx_l) \in \hat{H}^l$, and 
      integers $\vec{a} \in \BZ^l$, 
      such that  $\gx=\gx_k$ for some 
      $l \ge k \ge 1$ and 
      $(\gx_1 \circ \gr, \ldots ,\gx_l\circ \gr)$ is of 
      type $\vec{a}$ w.r.t . $\bar{\triangle}_{\vec{a}}(\mu_{Y})$.
      Let $Y$ be a $(j-1)$-step pro-nilflow, and let $(f_1,f_2, \ldots,f_{l})$
              be of type $\vec{a}$ w.r.t $\bar{\triangle}_{\vec{a}}(\mu_{Y})$.
       Then 
      \begin{enumerate}
        \item \label{item:cohomologus_nilflow} 
               $f_k$ is cohomologous to a function lifted from a 
               $(j-1)$-step nilflow.
        \item \label{item:countable} 
              For $k=1,\ldots,l$, $f_k$ belongs to a countable set modulo 
              coboundaries.
        \item \label{item:set_type_j_to_type_j} For $k=1,\ldots,l$, 
              $f_k:Y \ra S^1$ is of type $j$. 
        \item \label{item:type_j_to_type_j}
              If $\gr:Y \ra H$ is of type $j$ then for any character $\gx$ of 
              $H$, $\gx \circ \gr$ is of type $j$.
        \item \label{item:group} If $f,g:Y \ra S^1$ are of type $j$, then 
              $fg$ is of type $j$.
      \end{enumerate}
\item \label{item:type_j_pronil}
      If  $X=Y_{j}(X)\times_{\gs} H$ is an abelian extension by a cocycle
      of type $j$, then $X$ can be given the structure of a  
      $j$-step pro-nilflow. If $Y_j(X)$ is a nilflow and $H$ is a finite 
      dimensional torus, then  $X$ is a nilflow.  
\item \label{item:factor_nil}
      A factor of a $j$-step pro-nilflow is a $j$-step pro-nilflow.

\item \label{item:abelian_ext_type_j}
      If $X$ a $j$-step pro-nilflow then 
      $X=Y_{j}(X)\times_{\gs} H$ is an abelian extension of
      $Y_{j}(X)$ by a cocycle of type $j$. If $j\ge 1$ then $H$ is connected.
\item \label{item:main}
      $Y_{j+1}(X)$ can be given a structure of a $j$-step pro-nilflow.
\item \label{eq:convergence}  Let $a_1,\ldots a_{j+1} \in \BZ$, and 
      $f_1 \ldots f_{j+1} \in L^{\infty}(\mu_X)$. Then the averages 
       \begin{equation*}
            \Avr{N}{n} \p_{k=1}^{j+1} f_k(T^{a_kn}x)
       \end{equation*}
       converge in $L^2(\mu_X)$.
\end{enumerate}
\end{thm}

\begin{proof}
For $j=0$, $Y_{j+1}(X)$ is the pro-cyclic factor. For $j=1$, 
$Y_{j+1}(X)$ is the Kronecker factor which is an abelian extension
(by a connected group) of the pro-cyclic factor by a constant cocycle,
and all statements are easily verified.
Assume all statements hold replacing  $j$ with $j-1$.

\begin{dsc}\label{proof:def_mu_l}{\em proof of theorem \ref{thm:main} (\ref{def_mu_l})}

\[ 
\begin{split}
 \lim_{N \ra \infty} & \Avr{N}{n} \int \tau_{\vec{a}}^n f_1 \otimes \ldots \otimes
      f_{j+1} d\triangle_{j+1}(\mu_X)\\
 =& \lim_{N \ra \infty}\Avr{N}{n} \int f_1(T^{na_1}x) f_2(T^{na_2}x)\ldots 
    f_{j+1}(T^{na_{j+1}}x) d\mu_X^{}\\
 =& \lim_{N \ra \infty}\Avr{N}{n} \int f_1(x) f_2(T^{n(a_2-a_1)}x)\ldots 
    f_{j+1}(T^{n(a_{j+1}-a_1)}x) d\mu_X^{}\\
 =& \lim_{N \ra \infty} \Avr{N}{n} \int E(f_1|{Y}_{j})(\pi_{j}x) 
     \p_{i=1}^{j} (T^{n(a_{i+1}-a_1)}E(f_{i+1}|{Y}_{j})(\pi_{j}x)) d\mu^{}_{Y_{j}}\\
 =& \int E(f_1|Y_{j}) \otimes \ldots \otimes
          E(f_{j+1}|{Y}_{j})d\bar{\triangle}_{\vec{a}}(\mu_{Y_{j}}).
\end{split}
\]
By the above calculation
\begin{equation}\label{eq:conditional}
\int f_1 \otimes \ldots \otimes f_{j+1} d\bar{\triangle}_{\vec{a}}(\mu_{X}).=
\int E(f_1|{Y}_{j}) \otimes \ldots \otimes E(f_{j+1}|{Y}_{j})
     d\bar{\triangle}_{\vec{a}}(\mu_{Y_{j}}),
\end{equation}
thus 
$\bar{\triangle}_{\vec{a}}(\mu_{X})$ is the conditional product measure 
relative to $\bar{\triangle}_{\vec{a}}(\mu_{Y_{j}})$. 
\end{dsc}

\begin{dsc}\label{proof:iso}{\em proof of theorem \ref{thm:main}(\ref{item:iso})}.

We must show that if $E(f_k|\caB_{\hat{Y}_{j}(X)})=0$ for some $k$,
then the averages $\Avr{N}{n} \p_{k=1}^{j+1} f_k(T^{a_kn}x)$ converge to zero
in $L^2(\mu_X)$.
We apply the van der Corput Lemma \ref{VDC} with
\[
  u_n =\p_{k=1}^{j+1} T^{na_k} f_k(x).
\]
We calculate $\gc_m$:
\[
 \begin{aligned}
 \gc_m &= \lim \Avr{N}{n} \ip{u_n}{u_{n+m}}         \\
       &= \lim \Avr{N}{n} \int 
       \p_{k=1}^{j+1} T^{na_k}f_k(x) T^{na_k+ma_k}f_k(x)d\mu_X^{}\\
       &= \lim \Avr{N}{n} \int \p_{k=1}^{j+1} T^{na_k}(f_k T^{ma_k}f_k(x))d\mu_X^{} \\
       &= \int \left(f_1 \otimes \ldots \otimes f_{j+1} \right)
     \tau_{\vec{a}}^{m}\left(f_1 \otimes \ldots \otimes f_{j+1}\right) 
     d\bar{\triangle}_{\vec{a}}(\mu_X).
 \end{aligned}
\] 
By the ergodic theorem, there exists a $\tau_{\vec{a}}$ invariant function $D_{\vec{a}} \in  
L^2(\bar{\triangle}_{\vec{a}}(\mu_X))$ such that
\begin{equation}\label{D}
  \lim \Avr{M}{m} \gc_m = \int f_1 \otimes \ldots \otimes f_{j+1} 
                          D_{\vec{a}}(x_1,\ldots,x_{j+1}) d\bar{\triangle}_{\vec{a}}(\mu_X).
\end{equation}
By \ref{proof:def_mu_l}, 
$\bar{\triangle}_{\vec{a}}(\mu_{X})$ is 
the conditional product measure relative to 
$\bar{\triangle}_{\vec{a}}(\mu_{Y_{j}})$.
By theorem \ref{max_isometric}, $D_{\vec{a}}$  is measurable w.r.t 
$(\hat{Y}_{j}(X))^{j+1}$. If $E(f_k|\caB_{\hat{Y}_{j}(X)})=0 $ 
then the average (\ref{D}) is zero, and by VDC so is the original average.
\end{dsc}

\begin{dsc} {\em proof of theorem \ref{thm:main}(\ref{def:type_j})}.

This part is the bulk of the theorem, and the proof of its items will
be intertwined with the proof of the rest of the items in 
theorem \ref{thm:main}.
Let $Y$  a $(j-1)$-step pro-nilflow, with $j \ge 2$. By 
corollary \ref{cor:ucf_pronil}, $Y=Y_j(Y)$. By
the induction hypothesis in theorem 
\ref{thm:main}(\ref{item:abelian_ext_type_j}), we  
can identify $Y$ with a presentation as a tower of abelian extensions
$Y=H_1 \times_{\gs_1} \times H_2 \times \ldots \times_{\gs_{j-1}} H_j$
where $\gs_i$ is of type $i$, $H_i$ is connected for $i>1$, and 
$Y_i(Y)=H_1 \times_{\gs_1} \times H_2 \times \ldots \times_{\gs_{i-1}} H_i$. 
Specifically   
$Y=Y_{j-1}(Y) \times_{\gs_{j-1}} H_j$, where $H_j$ is
a connected compact abelian group, and $\gs_{j-1}$ is of type $j-1$.
Denote $Y_{j-1}=Y_{j-1}(Y)$.
Let $\pi_{j-1}:Y \ra Y_{j-1}$ be the projection. 
We identify $y \in Y$ with  $(\pi_{j-1}y,h) \in \pi_{j-1}Y \times H_j$.
Let $l$ be a positive integer.
To simplify the notation we now restrict ourselves to the special case
where $\vec{a}=(1,2,\ldots,l)$. The analysis is similar for any 
$\vec{a} \subset \BZ^l$. 

We write $\bar{\triangle}_l(\mu_Y)$ 
for the measure $\bar{\triangle}_{(1,\ldots,l)}(\mu_Y)$, and we say 
that $(f_1,\ldots,f_{l})$ is of type $l$
w.r.t $\bar{\triangle}_l(\mu_Y)$ if $(f_1,\ldots,f_{l})$ is of type 
$(1,\ldots,l)$ w.r.t $\bar{\triangle}_{l}(\mu_Y)$. In this case  
there exists a function $F \in L^{\infty}(\bar{\triangle}_l(\mu_Y))$ 
such that 
\begin{equation}\label{eq:func_j+1}
  \p_{k=1}^{l} f_k^{(k)}(y_k)
      =\frac{\tau F(y_1\ldots,y_{l})}
            {F(y_1\ldots,y_{l})}.
\end{equation}
\end{dsc}
\begin{rmr}\label{rmr:extra_coordinate}
As ${Y}$ is a $(j-1)$-step pro-nilflow,
on the support of $\mu_l(Y)$, the coordinates
$y_{j+1},\ldots,y_{l}$  are determined by the first 
$j$ coordinates $y_{1},\ldots,y_{j}$, and this correspondence 
is invariant under
$\tau$ (if $j=2$ then $Y$ is an abelian group, and 
$y_1,y_2,y_3$ form an arithmetic sequence, in general 
see the discussion in \ref{dsc:ergodic_com}). 
Therefore the function $F(y_1,\ldots,y_{l})$
can be replaced by a function of $j$ coordinates, and
equation (\ref{eq:func_j+1}) can be written in the form
\begin{equation}\label{eq:Func_j+1}
  \p_{k=1}^{l} f_k^{(k)}(y_k)
      =\frac{\tau  F(y_1,\ldots,y_{j})}
            {F(y_1,\ldots,y_{j})}.
\end{equation}
We will repeatedly refer to this equation.
\end{rmr}

\begin{rmr}\label{rmr:freedom_in_H_j}
The measure $\bar{\triangle}_{j}(\mu_Y)=
\bar{\triangle}_{j}(\mu_{Y_{j-1}}) \times (m_{H_j})^{j}$
(replace $X$ by $Y$ and $j$ by $j-1$ in equation (\ref{eq:conditional})).
\end{rmr}

\begin{lma}\label{lma:g_u_k} Let $(f_1,\ldots,f_{j+1})$ be of type $j+1$ w.r.t 
$\bar{\triangle}_{j+1}(\mu_Y)$. Then 
for each $k=1,\ldots,j+1$ there exists a family of functions 
$\{g_{k,u}\}_{u \in H_{j}} \subset B(Y^0,S^1)$, and 
a family of functions $\{f_{k,u} \}_{u \in H_{j}} \subset B(Y^0,S^1)$ 
such that 
\begin{equation}\label{eq:g_u_k}
  \frac{f_k(\pi_{j-1} y,h+u)}
       {f_k(\pi_{j-1} y,h)}
  = g_{k,u}(\pi_{j-1} y) \frac{Tf_{k,u}(y)}{f_{k,u}(y)}.
\end{equation}
\end{lma}

\begin{proof}
We use the fact that $\bar{\triangle}_{j+1}(Y)$ is 
invariant under translations by elements of the Mackey group $M_{j-1,j+1}$, 
and by elements of  
$\triangle_{j+1}(H_j)=\{(h,\ldots,h)\}_{h \in H_j} \subset H_j^{j+1})$.
This group is described in \ref{abelian_ext_nil}.
Let 
\[
M_{j-1,j+1}(0)=(M_{j-1,j+1}+\triangle_{j+1}(H_j)) \cap (H_j^{j} \times \{0\})
\]
The projection of $M_{j-1,j+1}(0)$ on any $j$ 
coordinates is full (i.e $H_j^j$).
Let $\vec{u}=(u_1,\ldots,u_{j+1}) \in M_{j-1,j+1}(0)$ ($u_{j+1}=0$). Then 
\begin{equation}\label{eq:Func_j+1_u}
  \p_{k=1}^{j+1} f_k^{(k)}(\pi_{j-1} y_k,h_k+u_k)
   =\frac{\tau F(\pi_{j-1} y_1,h_1+u_1,\ldots,\pi_{j-1} y_{j},h_{j}+u_{j})}
            {F(\pi_{j-1} y_1,h_1+u_1,\ldots,\pi_{j-1} y_{j},h_{j}+u_{j})}.
\end{equation}
Dividing equation (\ref{eq:Func_j+1_u}) by equation 
(\ref{eq:Func_j+1}) we get
\begin{equation}\label{eq:g_u_k_1}
  \p_{k=1}^{j} \frac{f_k^{(k)}(\pi_{j-1} y_k,h_k+u_k)}
                            {f_k^{(k)}(\pi_{j-1} y_k,h_k)}
   =\frac{\tau F_{\vec{u}}(y_1,\ldots,y_{j})}
            {F_{\vec{u}}(y_1,\ldots,y_{j})},
\end{equation}
where 
\[
F_{\vec{u}}(y_1,\ldots,y_{j})=
 \frac{F(\pi_{j-1} y_1,h_1+u_1,\ldots,\pi_{j-1} y_{j},h_{j}+u_{j})}
      {F(y_1,\ldots,y_{j+1})}
\]
Using remark \ref{rmr:freedom_in_H_j} we can apply  
theorem \ref{thm:hom_general} to $\bar{\triangle}_{j}(\mu_Y)$.
For any $1\le k \le j$ and any $u_k \in H_j$ there exist functions 
$g_{k,u_k} \in B(Y^0,S^1)$ and $f_{k,u_k} \in B(Y^0,S^1)$ such 
that 
\begin{equation}\label{eq:g_u_k_2}
  \frac{f_k^{(k)}(\pi_{j-1} y,h+u_k)}{f_k^{(k)}(\pi_{j-1} y,h)} 
= g_{k,u_k}(\pi_{j-1} y) \frac{T^k f_{k,u_k}(y)}{f_{k,u_k}(y)}.
\end{equation}
We still need to show that the same holds for $f_k$ (rather than $f_k^{(k)}$):
we use the fact that $T_{f_k}$ (see \ref{isometric_extensions}) 
and $T^k_{f_k^{(k)}}=(T_{f_k})^k$ commute
therefore 
\begin{equation}\label{eq:T_T^k_commute}
f_k^{(k)}(Ty)f_k(y)=f_k(T^ky)f_k^{(k)}(y).
\end{equation}
Using equations (\ref{eq:g_u_k_2}),(\ref{eq:T_T^k_commute}), and a calculation
we find that the function
\[
\frac{f_k(\pi_{j-1} y,h+u+v)/f_k(\pi_{j-1} y,h+u)}
           {f_k(\pi_{j-1} y,h+v)/f_k(\pi_{j-1} y,h)}\left/
\frac{T(f_{k,v}(\pi_{j-1} y,h+u)/f_{k,v}(\pi_{j-1} y,h))}
           {f_{k,v}(\pi_{j-1} y,h+u)/f_{k,v}(\pi_{j-1} y,h)} \right.
\]
is $T^k$ invariant and therefore constant  on the 
(finitely many) ergodic components of $T^k$. Denote this constant
$\gd_{k,u,v}(\pi_1y)$ (any ergodic component of $T^k$ is determined by a
cyclic group which is a factor of $Y_1(X)$ same argument as in 
\ref{thm:hom_general}).
By lemma \ref{multiplicity_of_f_h}, $\gd_{k,u,v}(\pi_1y)$ is multiplicative in 
$u$ (also in $v$) in a neighborhood of zero in $H_j$. 
By equation (\ref{eq:g_u_k_2}) (after iteration)
$(\gd_{k,u,v})^k(\pi_1y)$ is an eigenvalue, therefore $\gd_{k,u,v}(\pi_1y)=1$ 
for $u$  in a 
neighborhood of zero in $H_j$. 
Iterating we find this is true for all $u \in H_j$ ($H_j$ is connected). 
By corollary \ref{hom}
\[
\frac{f_k(\pi_{j-1} y,h+v)}{f_k(\pi_{j-1} y,h)}= \ti{g}_{k,v}(\pi_{j-1} y) 
\frac{T\ti{f}_{k,v}(\pi_{j-1} y,h)}{\ti{f}_{k,v}(\pi_{j-1} y,h)}.
\]
\end{proof}

\begin{lma}\label{lma:lift}
If $(g_1 \circ \pi_{j-1},\ldots, g_{l} \circ \pi_{j-1})$ is of type $l$ w.r.t
$\bar{\triangle}_{j}(\mu_{Y})$, then  as a function on 
$\pi_{j-1}Y$, $g_k$ is of type $j-1$, for $k=1,\ldots,l$.
\end{lma}

\begin{proof}
By definition
there exists a $\bar{\triangle}_{l}(\mu_{Y})$ measurable function $L$
such that
\[
  \p_{k=1}^{l} (g_{k})^{(k)}(\pi_{j-1} y_k)
      =\frac{\tau  L(y_1,\ldots,y_{l})}
            {L(y_1,\ldots,y_{l})}.
\]
Taking the Fourier expansion of $L$ with respect to the abelian group $H_j^{j}$
\[
L(y_1,\ldots,y_{l})
 =\sum_{\gx \in \hat{H}_j^{j}}G_{\gx}(\pi_{j-1} y_1,\ldots,\pi_{j-1} y_{l})
   \gx_1(h_1)\ldots \gx_{j}(h_{l}),
\]
We find that for any $\gx \in \hat{H}_j^{l}$
\[
G_{\gx}(y_1,\ldots,y_l)
\p_{k=1}^{l} g_{k}^{(k)}(\pi_{j-1} y_k)\bar{\gx}_{k}(\gs_{j-1}^{(k)}(\pi_{j-1} y_k))
=\tau G_{\gx}(\pi_{j-1} y_1,\ldots,\pi_{j-1} y_{l}).
\]
The function $|G_{\gx}|$ is invariant under $\tau$ and therefore
constant on the ergodic components of $\bar{\triangle}_{l}(\mu_{Y_{j-1}})$ . 
As $\hat{H}_j^{l}$
is countable, there is a character $\gx \in \hat{H}_j^{l}$,
and a set of $\bar{\triangle}_{l}(\mu_{Y_{j-1}})$ positive measure $A$, 
that is 
$\tau$ invariant, and for which
\begin{equation}\label{eq:A}
\p_{k=1}^{l} g_{k}^{(k)}(\pi_{j-1} y_k)\bar{\gx}_{k}
  \gs_{j-1}^{(k)}(\pi_{j-1} y_k)
=\frac{\tau G_{\gx}(\pi_{j-1} y_1,\ldots,\pi_{j-1} y_{l})}
      {G_{\gx}(\pi_{j-1} y_1,\ldots,\pi_{j-1} y_{l})}.
\end{equation}

Denote  by $W$ the system 
$(Y_{j-1}^l,\bar{\triangle}_{l}(\mu_{Y_{j-1}}),\tau)$. 
Let $\gr_k=g_{k}\bar{\gx}_{k}(\gs_{j-1})
:Y_{j-1} \ra S^1$, and let  
$\ti{\gr}=(\gr_1,\gr^{(2)}_2, \ldots, \gr^{(l)}_{l})$. Denote $\bar{y}
:=\pi_{j-1}y$.
Consider the group extension 
$
W \times_{\ti{\gr}} (S^1)^{l}.
$
Let $W_{\bar{y}}$ be an ergodic component of $W$ 
(the ergodic components of $W$ are parametrized by $\bar{y} \in Y_{j-1}$), and 
let $P_{\bar{y}} \subset (S^1)^{l}$ be the Mackey group associated with 
the ergodic components of the group extension 
$W_{\bar{y}} \times_{\ti{\gr}} (S^1)^{l}$.
As the transformations 
\[ 
T_{\gr_1}\times (T_{\gr_2})^{2}
\times \ldots \times (T_{\gr_l})^{l}, \quad
T_{\gr_1}\times T_{\gr_2}
\times \ldots \times T_{\gr_l}
\]
commute, $P_{T\bar{y}}=P_{\bar{y}}$ (see lemma \ref{lma:conjugate}), 
and by ergodicity
$P_{\bar{y}}=P$ is constant  a.e. as a function of $\bar{y}$. 
As equation (\ref{eq:A}) holds on a $\tau$ invariant set of 
$\bar{\triangle}_{l}(\mu_{Y_{j-1}})$ positive measure, 
we get
\[
(\gf,\ldots,\gf) \in P^{\perp},
\]
where $\gf(\gz)=\gz$. 
Therefore for a.e. $\bar{y}$, there exists a function $G_{\bar{y}}$, so that 
equation (\ref{eq:A}) holds on $W_{\bar{y}}$ replacing $G_{\gx}$ with 
$G_{\bar{y}}$. 
Notice that  $G_{\bar{y}}$ is determined up to a constant multiple on 
$W_{\bar{y}}$, and 
can therefore be chosen so that it depends measurably on $\bar{y}$ 
(see \ref{g_x_constant}).
Thus there exists a measurable function $G$ such that 
equation (\ref{eq:A}) holds $\bar{\triangle}_{l}(\mu_{Y_{j-1}})$ a.e.
replacing $G_{\gx}$ by $G$.
This implies that  
$(g_{1} \cdot \bar{\gx}_{1}(\gs_{j-1}), \ldots ,
  g_{l} \cdot \bar{\gx}_{l}(\gs_{j-1}))$
is of type $l$ w.r.t. $\bar{\triangle}_{l}(\mu_{Y_{j-1}})$.
By the induction hypothesis in 
\ref{thm:main}(\ref{item:set_type_j_to_type_j}),
(\ref{item:type_j_to_type_j}) the functions 
$g_{k} \gx_{k}(\gs_{j-1})$, $\gx_{k}(\gs_{j-1})$ are of type $j-1$.
By the induction hypothesis in \ref{thm:main}(\ref{item:group}), 
$g_{k}$ is of type $j-1$.

\end{proof}

\begin{cor}\label{cor:lift}
Let $(f_1,\ldots,f_{j+1})$ be of type $j+1$
w.r.t $\bar{\triangle}_{j+1}(\mu_Y)$. 
Let  $\vec{u} \in M_{j-1,j+1}(0)$, and let $g_{k,u_k}$ satisfy
equation (\ref{eq:g_u_k}). Then for $k=1,\ldots j+1$, as a function on 
$\pi_{j-1}Y$, $g_{k,u_k}$ is of type $j-1$. 
\end{cor}
\begin{proof} Substitute equation 
(\ref{eq:g_u_k_2}) in  equation (\ref{eq:g_u_k_1}) and use lemma  
\ref{lma:lift}.
\end{proof}

\begin{cor}\label{cor:center} If $(f_1,\ldots f_{j+1})$ is of type $j+1$
w.r.t $\bar{\triangle}_{j+1}(\mu_Y)$ then for $k=1,\ldots ,j+1$
there exists a family of constants 
$\{\gl_{k,u}\}_{u \in H_{j}}$, and 
a family of measurable functions $\{f_{k,u} \}_{u \in H_{j}}$ such that 
\begin{equation}\label{center_eq}
  \frac{f_k(\pi_{j-1} y,h+u)}
       {f_k(\pi_{j-1} y,h)}
  = \gl_{k,u} \frac{Tf_{k,u}(y)}{f_{k,u}(y)}.
\end{equation}
\end{cor}

\begin{proof}
By theorem \ref{thm:main}(\ref{item:countable}), using the induction hypothesis, there are countably many $g_{k,u}$ up to $\pi_{j-1}Y$-quasi-coboundaries. 
There exists a set $U$ of 
positive measure in $H_j$ such that
\[
 u,v \in U \Rightarrow \frac{g_{k,u}}{g_{k,v}} \ \te{is a quasi-coboundary}.
\] 
If $u,u+v \in U$ and  if $f_{k,u,v}=f_{k,u+v}/f_{k,u}$ then 
\[
 \begin{split}
\frac{f_k(\pi_{j-1} y,h+u+v)}{f_k(\pi_{j-1} y,h+u)}  
= \frac{f_k(\pi_{j-1} y,h+u+v) / f_k(\pi_{j-1} y,h)}
       {f_k(\pi_{j-1} y,h+u) / f_k(\pi_{j-1} y,h)}
= C_{u,v,k} \frac{T f_{k,u,v}(y)}{f_{k,u,v}(y)}.
\end{split}
\]
Thus the claim is true for $v$ in a neighborhood of zero in $H_j$
(as the map $H \ra H$ sending $h$ to $h+u$ is onto). As $H_j$
is connected equation (\ref{center_eq}) holds for all $v \in H_j$.
\end{proof}

\begin{lma}\label{lma:mult_of_f_u}
 The families in the previous lemma can be chosen so that 
the function $\gl_{k,u}:H_{j} \ra S^1$ is multiplicative in a neighborhood of 
zero in $H_{j}$.
\end{lma}

\begin{proof}
By lemma \ref{multiplicity_of_f_h}.
\end{proof}

\begin{cor}\label{cor:product_lamda}
The functions $f_{k,u}$  can be chosen so that for some neighborhood of zero
$U \subset H$,  for any 
$\vec{u}=(u_1, \ldots ,u_{j+1}) \in M_{j-1,j+1}+\triangle_{j+1}(H_j) \cap 
U^{j+1}$ 
\[\p_{k=1}^{j+1} \gl^k_{k,u_k}=1\]
\end{cor}

\begin{proof}
Choose the families $\{f_{k,u}\},\{\gl_{k,u}\}$ so that $\gl_{k,u}$
is multiplicative in a a neighborhood of 
zero in $H_{j}$. 
Substituting equation (\ref{center_eq}) in equation 
(\ref{eq:g_u_k_1}), we find that $\p_{k=1}^{j+1} \gl^k_{k,u_k}$ is an 
eigenvalue of $\tau$.
\end{proof}

\begin{lma}\label{lma:H_j_finite}
Let $(f_1,\ldots,f_{j+1})$ be of type $j+1$ w.r.t 
$\bar{\triangle}_{j+1}(\mu_Y)$. Then there exists an integer $n$, 
and a factor $\ti{Y}=Y_{j-1}\times \BT^n$ of $Y$, such that if 
$p:Y\ra \ti{Y}$ is the factor map, then for $k=1,\ldots,j+1$,
$f_k$ is cohomologous to a cocycle $\ti{f}_k \circ  p$. 
Furthermore, there exist functions
$g_1,\ldots,g_{j+1}$, where $g_k:Y^0_{j-1} \ra S^1$ is of type $j-1$,
such that 
\[
(\ti{f}_1g_1 \circ \pi_{j-1},\ldots,
\ti{f}_{j+1}g_{j+1}\circ \pi_{j-1})
\] 
is of type $j+1$ w.r.t
$\bar{\triangle}_{j+1}(\mu_{\ti{Y}})$, and therefore corollary 
\ref{cor:center} holds replacing $f_k$ by $\ti{f}_k$, and $Y$ by $\ti{Y}$
($\pi_{j-1}$ is the projection $\ti{Y} \ra Y_{j-1}$).
\end{lma}

\begin{proof}
By lemma \ref{reduction_finite_connected}, and remark \ref{reduction_no_torus}
we can find $J$  so that
$H_{j}/J=\BT^n \times \ti{H}$ where $\ti{H}$ is a compact totally disconnected 
abelian group (lemma \ref{reduction_finite_connected} can be carried out 
simultaneously for $(f_1,\ldots,f_{j+1})$). Consider the system 
$W=Y_{j-1} \times_{\gs_{j-1}} (\BT^n \times \ti{H})$. 
This system is a factor of $Y_{j}$, therefore  $Y_{j-1}(W)=Y_{j-1}$, 
$W$ is a $(j-1)$-step pro-nilflow, and  $Y_{j}(W)=W$ 
(this follows from the induction hypothesis in 
theorem \ref{thm:main}(\ref{item:factor_nil}), and from corollary \ref{cor:ucf_pronil}).
This  a contradiction
to $Y_{j}(W)$ being an extension of $Y_{j-1}(W)$ by a connected abelian group
for $j>1$ (this follows from the induction hypothesis in 
\ref{thm:main}(\ref{item:abelian_ext_type_j})). Now $Y$ can be presented as a 
skew product $Y=\ti{Y} \times_{\ti{\gs}_{j-1}} J$ 
where 
\[
\ti{\gs}_{j-1}(\ti{y})=\bar{r}(\ti{y})\gs_{j-1}(\pi_{j-1}(\ti{y}))r(T\ti{y}),
\] 
where $r(\ti{y})$ takes values in $H_j$. As in the proof of lemma 
\ref{lma:lift}, 
there exists a character $(\gx_1,\ldots,\gx_{j+1}) \in \hat{J}^{j+1}$ such 
that 
\[
(\ti{f}_{1} \cdot {\gx}_{1}(\ti{\gs}_{j-1}), \ldots ,
  \ti{f}_{j+1} \cdot {\gx}_{j+1}(\ti{\gs}_{j-1}))
\]
is of type $j+1$ w.r.t. $\bar{\triangle}_{j+1}(\mu_{\ti{Y}})$.
Let $\ti{\gx}_k$ be the lift of $\gx_k$ to a character of $H_j$,
then 
\[
  \ti{\gx}_k(\ti{\gs}_{j-1}(\ti{y}))
   =\ti{\gx}_k (\bar{r}(\ti{y})) \ti{\gx}_k(\gs_{j-1}(\pi_{j-1}(\ti{y})))
       \ti{\gx}_k (r(T\ti{y}).
\]
Therefore 
\[
(\ti{f}_{1} \cdot {\ti{\gx}}_{1}(\gs_{j-1}(\pi_{j-1})), \ldots ,
  \ti{f}_{j+1} \cdot {\ti{\gx}}_{j+1}(\gs_{j-1}(\pi_{j-1}))
\]
is of type $j+1$ w.r.t. $\bar{\triangle}_{j+1}(\mu_{\ti{Y}})$.
\end{proof}

\begin{ntt}
If $U$ is an abelian group we denote rotation by an element of $U$
by $R_u$.
\end{ntt}

\begin{lma}\label{finite_dim_nil} Let $(f_1,\ldots,f_{j+1})$ be of type $j+1$ 
w.r.t 
$\bar{\triangle}_{j+1}(\mu_Y)$. Then  there exists a 
factor of $Y$ which is a $(j-1)$-step nilflow $\ti{Y}=(N/\gC,a)$, such that 
if $p:Y \ra \ti{Y}$ is the factor map, then for $k=1,\ldots j+1$
\begin{enumerate}
\item $f_k$ is cohomologous to $\ti{f}_k \circ p$. 
\item There exist functions
$g_1,\ldots,g_{j+1}$, where $g_k:\ti{Y}^0 \ra S^1$ is of type $j-1$,
such that 
$
(\ti{f}_1g_1 \circ \pi_{j-1},\ldots,
\ti{f}_{j+1}g_{j+1}\circ \pi_{j-1})$ 
is of type $j+1$ w.r.t
$\bar{\triangle}_{j+1}(\mu_{\ti{Y}})$
($\pi_{j-1}$ is the projection $\ti{Y} \ra \pi_{j-1}(\ti{Y})$)
\item Corollary \ref{cor:center} holds replacing 
     $f_k$ by $\ti{f}_k$, and $Y$ by $\ti{Y}$.
\end{enumerate}
\end{lma}

\begin{proof} Recall the identification of  
$Y$ as a tower of abelian extensions
$Y=H_1 \times_{\gs_1} H_2 \times \ldots \times_{\gs_{j-1}} H_j$,
where $Y_i(Y)=Y_{i-1}(Y) \times_{\gs_{i-1}} H_i$,
the group $H_i$ is an abelian 
group which is connected for $i>1$, and $\gs_{i}$ is of type $i$ .
We would like to 'replace' $H_i$ for $i\ge 2$ with a finite dimensional 
torus $\BT^{n_i}$ (and replace $H_1$ by a cyclic group) and get a system
$C \times_{\gs'_{1}} \BT^{n_2}\times_{\gs'_{2}} \ldots \times_{\gs'_{j-1}} \BT^{n_j}$ that is a factor of $Y$ and therefore a nilflow 
(a pro-nilflow by the induction hypothesis in theorem \ref{thm:main}(\ref{item:factor_nil}), and a nilflow by the induction hypothesis in 
\ref{thm:main}(\ref{item:type_j_pronil})), 
and if $p$ is the factor map, then  $f_k$ is cohomologous to 
$\ti{f}_k \circ p$. 

We do this by decreasing 
induction on the index $i$. The case $i=j$ 
was proved in lemma \ref{lma:H_j_finite}. 
Assume we have constructed a system 
$\ti{Y}=H_1 \times_{\gs_{1}} \ldots \times H_i\times_{\gs'_i} \BT^{n_{i+1}}
\times \ldots \times_{\gs'_{j-1}} \BT^{n_j}$ that is a factor of $Y$
and satisfies $(1)-(3)$.  We may now forget the original system $Y$. By 
abuse of notation we replace $\ti{Y}$ by $Y$, and $\ti{f}_k$ by $f_k$.  
The cocycles 
$\gs'_i, \ldots,\gs'_{j-1}, g_1,\ldots g_{j+1}$ are of type $< j$ and 
take values in finite 
dimensional tori. Therefore there exists a finite dimensional torus 
$\BT^{n_i}$ with $\gs'_l$ for $l=i,\ldots ,j-1$, 
cohomologous to 
functions $\ti{\gs_l}$,  lifted from  
$Y'_l:=H_1 \times\ldots \times H_{i-1}\times_{\gs'_{i-1}} \BT^{n_{i}}
           \times \ldots \times_{\gs'_{l-1}} \BT^{n_l}$,
and $g_k$, for $k=1,\ldots j+1$ cohomologous to $\ti{g}_k$ 
lifted from $Y':=Y'_{j+1}$. 
After reparametrization we may assume $\gs'_l$, for $l=i,\ldots ,j-1$,
is lifted from $Y'_l$. Condition $(2)$ remains valid if we can replace 
$g_k$ by $\ti{g}_k$.  
Write $H_i =\BT^{n_i}\times U$ (measurewise) where $U$
is a compact abelian group.
The action of $T_Y$ 
commutes with
rotation by an element in $U$; i.e with
\[
R_{u'}:(h_1,\ldots,h_{i-1},t_i,u,t_{i+1},\ldots, t_j) \ra
(h_1,\ldots,h_{i-1},t_i,u+u',t_{i+1},\ldots, t_j)\] 
Indeed, 
$\gs'_i,\ldots, \gs'_{j-1}$ are not affected by translation in elements of 
$U$. Both $R_u, T_Y$
commute with rotation by an element $t \in \BT^{n_j}$.
Therefore
\[\frac{f_k(R_u (\pi_{j-1} y,t_j+t))}{f_k(R_u(\pi_{j-1} y,h))}\left/
  \frac{f_k(\pi_{j-1} y,t_j+t)}{f_k(\pi_{j-1} y,t_j)}\right.
= \frac{f_{k,v}(R_u Ty)}{f_{k,t}(Ty)}\left/
  \frac{f_{k,t}(R_uy)}{f_{k,t}(y)}\right. .
\]
By Theorem \ref{sthmCLhom} 
\[f_k(R_u y)/f_k(y)
\]
is cohomologous to a cocycle lifted from $\pi_{j-1}Y$. 
By the same argument as in lemma \ref{cor:center} 
it is cohomologous to a constant for $u$ in a neighborhood of zero in $U_i$.

Now proceed as in lemma \ref{lma:H_j_finite} to obtain $\ti{Y}$ 
(which will be a factor between $Y'$ and $Y$).
Applying the same procedure 
for $i=1$ gives the cyclic part. 
\end{proof}

\begin{rmr}\label{rmr:countable_nil}  Let $f_1,\ldots,f_{j+1}$ be as in
 lemma \ref{finite_dim_nil}. 
By remark \ref{rmr:co-tori} there are countably many possibilities for the 
groups $U_i$ and therefore countably many possibilities for the nilflow
$N/\gC$ (up to isomorphism).
\end{rmr}

\begin{dsc}\label{dsc:nil} 
Let  $(N/\gC,a)$ be the  $(j-1)$-step nilflow from lemma \ref{finite_dim_nil},
and let $p:Y \ra N/\gC$ be the projection.
Let $N_1=N$, and $N_{l+1}=[N,N_l]$ ($N_{j}=\{1\}$). 
We will show that if  the system 
$(N/\gC,a)$ has no finite non trivial factors, then
$N/\gC \times_{\ti{f}_k} S_1$ can be given the structure of a 
$j$-step nilflow. If the nilflow $(N/\gC,a)$  has  a 
finite factor $C_i$, then $C_i$
is an abelian group of order $i$ for some integer 
$i$. The nilflow $(N/\gC,a^i)$ has finitely many ergodic components, and 
rotation by $a$ induces an isomorphism between them. Each
ergodic component $X$ (with the action of $a^i$) is a $j-1$ step nilflow
with no nontrivial finite factors.
We will show that  the system $X \times_{\ti{f}_k} S^1$ isomorphic to a 
$j$-step nilflow $(M/\gL,b)$.
The system $N/\gC \times_{\ti{f}_k} S^1$ will then be isomorphic to a union
of $i$ isomorphic $j$-step nilflows $C_i \times M/\gL$, with 
the action of $T$ given by: for $i-1> k \ge 0: (k,m\gL) \ra (k+1,m\gL)$
and for $k=i-1: (i-1,m\gl) \ra (0,bm\gL)$. The group generated
by $\{(0,m), T\}_{m \in M}$ is $j$-step nilpotent, 
and acts transitively on $C_i \times M/\gL$.

Assume $(N/\gC,a)$ has no finite non trivial factors.
The $r$-u.c.f of $N/\gC$ is $Y_r(N/\gC)= N/N_r\gC$ (this follows from 
lemma \ref{Y_j_of_nil}).
The system $(N/\gC,a)$ can be presented as an abelian extension of a 
$(j-2)$-step
nilflow i.e. $(N/N_{j-1}\gC \times N_{j-1}/( N_{j-1}\cap \gC)$, and
we may assume that $N$ is simply connected. 
The group  $N_{j-1}$ is abelian, connected and 
simply connected (\cite{L02}), therefore isomorphic 
to $\BR^m$ for some $m$. Let  the action of $T$ on 
$N/N_2 \cong \BZ \times \BR^n$ be given by translation by $\ga$.
Then equation (\ref{center_eq})
becomes: 
for any $c \in N_{j-1}$, $k=1, \ldots, j+1$, $y \in N/\gC$ 
\begin{equation}
\frac{\ti{f}_k(cy )}{\ti{f}_k( y)}=\gl_{k,c} 
\frac{Tf_{k,c}(y)}{f_{k,c}(y)}.
\end{equation}
\end{dsc}

\begin{lma}\label{multiplicity_f_s} 
Let $(N/\gC,a)$ be the $(j-1)$-step nilflow
from \ref{dsc:nil}.

Let $f_1,\ldots, f_k$ be functions in  $B(N/\gC,S^1)$.
Let $\{\gl_{k,c}\}$ be a family of constants, $\{f_{k,c}\}_{c \in N_{j-1}}$ be a family of functions in  $B(N/\gC,S^1)$ such that 
\begin{equation}\label{eq:center_N}
\frac{f_k(cy)}{f_k(y)}=\gl_{k,c} 
\frac{f_{k,c}(ay)}{f_{k,c}(y)}
\end{equation}
for all $c \in N_{j-1}$. Then   
We can choose $f_{k,c},\gl_{k,c}$ in equation (\ref{eq:center_N})
such that 
\[
  f_{k,c_1}(k,c_2y) f_{k,c_2}(y) \sim f_{k,c_1c_2}(y).
\]
\end{lma}

\begin{proof} It follows from lemma \ref{multiplicity_f_s_nil}.

\end{proof}

\begin{lma}\label{matrix}

Let $f_k, f_{k,c}, \gl_{k,c}$ be from lemma \ref{multiplicity_f_s}.
Then for each $k=1,\ldots j+1$ 
there exists an integer matrix $L_k$,  
a neighborhood of zero  $U \subset  N_{j-1}$, and a family of functions 
$\{f_{k,c}\}_{c \in U}$, such 
that for all $c \in U$ 
\[
\frac{f_k(cy)}{f_k(y)} =e^{2 \pi i \ip{L_k\ga}{c}}\frac{f_{k,c}(ay)}{f_{k,c}(y)} .
\]
\end{lma}

\begin{proof}
It follows from lemma \ref{matrix_nil}. 
\end{proof}

\begin{cor}\label{cor:countable} Let $Y$ be a $(j-1)$-step pro-nilflow,
$(f_1,\ldots,f_{j+1})$ of type $j+1$ w.r.t. $\bt_{j+1}(\mu_Y)$.
Then modulo quasi-coboundaries, 
$f_k$ belongs to a countable set. 
\end{cor}

\begin{proof} By lemma \ref{finite_dim_nil}, $f_k$ is cohomologous to a 
function $\ti{f_k} $ lifted from
a nilflow, and by remark \ref{rmr:countable_nil} there are countably many
possibilities for this nilflow. Fix the nilflow. 
If  $(\ti{f}_1,\ldots,\ti{f}_{j+1})$,  $(\ti{f}'_1,\ldots,\ti{f}'_{j+1})$ 
share the integer 
matrices $L_1,\ldots,L_{j+1}$ from lemma \ref{matrix}, then 
by corollary \ref{hom}, $\ti{f}_{k}/\ti{f}'_{k}$ is cohomologous to a 
function on $\pi_{j-1}(N/\gC)$.
By lemma \ref{lma:lift}, as a function on $\pi_{j-1}Y$, $\ti{f}_k/\ti{f}'_k$
is of type $j-1$. 
By the induction hypothesis in theorem 
\ref{thm:main}(\ref{item:countable}), $\ti{f}_k/\ti{f}'_k$
belongs to a countable set modulo $\pi_{j-1}Y$  quasi-coboundaries.
(If we use condition $(2)$ in lemma \ref{finite_dim_nil} then we can get
that  $\ti{f}_k/\ti{f}'_k$
belongs to a countable set modulo $\pi_{j-1}(N/\gC)$  quasi-coboundaries).
\end{proof}

\begin{lma}\label{lma:g_u_k_l}  Let $Y$ be a $(j-1)$-step pro-nilflow, 
and let $(f_1,\ldots,f_{l})$ be of type $l$
w.r.t $\bt_{l}(\mu_Y)$. Then 
for each $k=1,\ldots,l$ there exists a family of constants
$\{\gl_{k,u}\}_{u \in H_{j}}$, and 
a family of functions $\{f_{k,u} \}_{u \in H_{j}} \subset B(Y^0,S^1)$ 
such that 
\begin{equation}\label{eq:g_u_k_l}
  \frac{f_k(\pi_{j-1} y,h+u)}
       {f_k(\pi_{j-1} y,h)}
  = \gl_{k,u}\frac{Tf_{k,u}(y)}{f_{k,u}(y)}.
\end{equation}
\end{lma}

\begin{proof} We use induction on $l$. 
The proof for $l\le j+1$ is given in corollary \ref{cor:center}.
Assume the statement holds for $l$; we show it for $l+1$. Let
$(f_1,\ldots,f_{l+1})$ be of type $l+1$ w.r.t $\bt_{l+1}(\mu_Y)$, and let
$M_{j-1,l+1}(0)=(M_{j-1,l+1}+\triangle_{l+1}(H_j))\cap H_j^l \times \{0\}$.
Let $\vec{u}=(u_1,\ldots,u_{l+1}) \in M_{j-1,l+1}(0)$ ($u_{l+1}=0$).
Then 
\[
\left( \frac{f_1(\pi_{j-1} y,h+u_1)}{f_1(\pi_{j-1} y,h)},\ldots,
  \frac{f_l(\pi_{j-1} y,h+u_l)}{f_l(\pi_{j-1} y,h)} \right)
\]
is of type $l$ and by the induction hypothesis
\[
\frac{f_k(\pi_{j-1} y,h+u_k+u) / f_k(\pi_{j-1} y,h+u)}
     {f_k(\pi_{j-1} y,h+u_k) / f_k(\pi_{j-1} y,h)}
=\gl_{k,u,u_k}\frac{Tf_{k,u,u_k}(y)}{f_{k,u,u_k}(y)}.
\]
The projection of  $M_{j-1,l+1}(0)$ on any coordinate 
$k\le l$ is full. By lemma \ref{lma:mult_of_f_u}, fixing $u$,
$\gl_{k,u,u_k}$ is multiplicative in $u_k$ in a neighborhood
of zero in $H_j$. By lemma \ref{matrix}, $\gl_{k,u,u_k}$
is determined by an integer matrix. The same holds interchanging the roles of
$u$ and $u_k$. Therefore 
$\gl_{k,u,u_k} \equiv 1$ in a neighborhood of zero in $H_j$.
By lemma \ref{hom}
\[
\frac{f_k(\pi_{j-1} y,h+u)}{ f_k(\pi_{j-1} y,h)}=
 g_{k,u}(\pi_{j-1} y)\frac{Tf_{k,u}(y)}{f_{k,u}(y)}.
\] 
As in lemma \ref{lma:lift} and corollary \ref{cor:center}, the functions 
$g_{k,u}(\pi_{j-1} y),f_{k,u}(y) $
can be chosen so that $g_{k,u}(\pi_{j-1} y)$ is a constant function on 
$Y_{j-1}$.
\end{proof}

\begin{cor}\label{cor:any_l} We may replace the index $j+1$  in 
lemmas/corollaries 
\ref{lma:g_u_k} - \ref{cor:countable} by the index $l$ for any $l \ge j+1$.
\end{cor}

As a corollary we get
\begin{dsc}{\em Proof of theorem \ref{main}(\ref{item:cohomologus_nilflow})}:
It follows from lemma \ref{finite_dim_nil} and corollary \ref{cor:any_l}.
\end{dsc}

\begin{dsc} {\em Proof of theorem \ref{main}(\ref{item:countable})}:
It follows from corollaries \ref{cor:countable} and \ref{cor:any_l}.
\end{dsc}

\begin{dsc}
We now fix $N/\gC$ - a $(j-1)$-step nilflow as in \ref{dsc:nil}. Let 
$(f_1,\ldots,f_{l})$ be of  type $l$ w.r.t $\bt_{l}(\mu_{N/\gC})$.
We will show that the system $N/\gC \times_{f_k} S^1$ is isomorphic to a $j$ 
step nilflow. We have constructed families of functions
$\{f_{k,c}\}_{c \in N_{j-1}}$ 
satisfying equation (\ref{eq:center_N}).
\end{dsc}

\begin{lma}\label{commuting_center}
For any $c_1,c_2 \in N_{j-1}$, 
$f_{k,c_1}$ and $f_{k,c_2}$ 
satisfying equation (\ref{eq:center_N}) we have
\[
 \frac{f_{k,c_1}(c_2y)}{f_{k,c_1}(y)}=\frac{f_{k,c_2}(c_1y)}{f_{k,c_2}(y)}
\] 
\end{lma}

\begin{proof} It follows from lemma \ref{commuting_center_nil}.

\end{proof}

\begin{dsc} 
Consider the group
\[
\caG=\{(n,f): \ n \in N, \ f \in B(N/\gC,S^1)\},
\]
with  multiplication
\[
  (n,f)(m,g)=(nm,f^m g), \qquad  (f^m g)(y)=f(my)g(y).
\]
The elements of the form $(1,C)$ where $C$ is a constant are in $Z(\caG)$ - the center
of $\caG$.
For $c \in N_{j-1}$, we can interpret equation (\ref{eq:center_N}) as 
\[ [(a,f),(c,f_{c})]=(1,\gl_{c}) \in Z(\caG)
\]
We want to think of $f_{c}$ as 
elements of the $(j-1)$th  subgroup in the upper central series of $\caG$, 
which will be a 
$j$-step nilpotent group.
We now follow the 
derived series upward 
to construct for each element in $N$ a function $f_{k,n}$ that will satisfy 
good commuting relations with $(a,f_k)$.
\end{dsc}

\begin{ntt} For $n \in  N_i \backslash N_{i+1}$, $|n|=i$, and 
$f^n(y)=f(ny)$
\end{ntt}

\begin{pro}\label{P:construction_f_n}
Let $Y=(N/\gC,a)$ be a $(j-1)$-step nilflow. Let $(f_1,\ldots,f_{l})$
be of type $j$ w.r.t $\bt_l(\mu_Y)$.
For $k=1,\ldots,l$  there exists a family of functions  
$\caF_k=\{f_{k,n}\}_{n \in N}$,
where $f_{k,n}:N/\gC \ra S^1$ are measurable functions satisfying the following
conditions:
\begin{enumerate}
\item\label{I:commutation_a_n}
 \[
  \frac{f_k(ny)}{f_k(y)}\sim f_{k,[a,n]}(nay)\frac{f_{k,n}(ay)}{f_{k,n}(y)},
 \]
\item\label{I:commutation_c}
 For any  $c \in N_{j-1}$, $|n|>1$ : 
 \[
  f_{k,c}^n f_{k,n}= f_{k,n}^c f_{k,c}.
 \]
\item\label{I:commutation_f_n} 
 \[
   \frac{f_{k,n_1}(n_2y)}{f_{k,n_1}(y)} \sim f_{k,[n_1,n_2]}(n_2n_1y)
   \frac{f_{k,n_2}(n_1y)}{f_{k,n_2}(y)},
 \]
\item\label{I:mult_n_1_n_2}
 \[
  f_{k,n_1n_2}(y) \sim  f_{k,n_1}(n_2y)f_{k,n_2}(y).
 \]
\item\label{I:invariance}
 If $(n_1,\ldots,n_{l}) \subset N^{l}$ 
 preserves the ergodic components of $\tau_{l}(Y)$ then 
 \begin{equation*}
  \frac{F(n_1y_1,\ldots,n_{l}y_{l})}
       {F(y_1,\ldots,y_{l})}
  \p_{k=1}^{l}\bar{f}_{k,n_k}(y_{k})
 \end{equation*}
is constant $\bt_{l}(\mu_Y)$ a.e.  
\end{enumerate}
\end{pro}

\begin{proof}
The proof requires a series of lemmas and their corollaries and will be completed in \ref{f_n_functional}.
We prove this inductively, proceeding upward in the derived series of 
$N$. By lemmas \ref{matrix}, 
\ref{multiplicity_f_s}, \ref{commuting_center},  and corollary 
\ref{cor:product_lamda},
conditions (\ref{I:commutation_a_n})-(\ref{I:mult_n_1_n_2})
 hold for all $n \in N_{j-1}$. 
As $N_{j-1} \subset Z(N)$ the function 
in (\ref{I:invariance}) is invariant under 
$\tau_{l},T_{l}$.
Suppose we constructed $f_n$ for $n$ in $N_{i+1}$ ($i+1>1$) satisfying 
conditions  (\ref{I:commutation_a_n})-(\ref{I:invariance}). 
Let $n \in N_i$.
Using conditions (\ref{I:commutation_a_n}),(\ref{I:commutation_c}),
for $c \in N_{j-1}$ 
we have:
\[
 \frac{f_k(ncy)}{f_k(cy)}\bar{f}_{k,[a,n]}(nacy)
      \left/\frac{f_k(ny)}{f_k(y)}\bar{f}_{k,[a,n]}(nay)\right.
 =\frac{f_{k,c}(nay)}{f_{k,c}(ay)}\left/\frac{f_{k,c}(ny)}{f_{k,c}(y)}\right. .
\]
By corollary \ref{hom} there exist functions $f_{k,n}:N/\gC \ra S^1$, and
$g_{k,n}:N/(N_{j-1}\gC)$  $\ra S^1$ 
such that 
\begin{equation}\label{g_n_still_there}
      \frac{f_k(ny)}{f_k(y)}\bar{f}_{k,[a,n]}(nay)
      =g_{k,n}(\pi_{j-1} y)\frac{f_{k,n}(ay)}{f_{k,n}(y)}.
\end{equation}
Let $c \in N_{j-1}$. By the induction hypothesis 
$f_{k,c}^{[a,n]}f_{[a,n]}=f_{[a,n]}^cf_{k,c}$ therefore
\[
\frac{f_{k,c}(ny)f_{k,n}(y)}{f_{k,c}(y)f_{k,n}(cy)}
\]
is a $T$ invariant function and therefore constant $\gd(k,n,c)$.

\begin{lma}\label{g_l_type_l} The $l$-tuple
$(g_{1,n},\ldots,g_{l,n})$ is of type $l$ w.r.t 
$\bt_l(\mu_Y)$.
\end{lma}

\begin{proof}
Iterating equation (\ref{g_n_still_there}) 
(using condition (\ref{I:mult_n_1_n_2})) we get in a neighborhood of zero in 
$N_i$
\[
      \frac{f_k^{(k)}(ny)}{f_k^{(k)}(y)}\bar{f}_{k,[a^k,n]}(na^ky)
      =g_{k,n}^{(k)}(\pi_{j-1} y)\frac{f_{k,n}(a^ky)}{f_{k,n}(y)}.
\]
Substituting in the functional equation
(\ref{eq:mackey}) we get
\begin{equation*} 
 \begin{split}
  \p_{k=1}^{l} & \left((g_{k,n}^{(k)}(\pi_{j-1} y_k)
   f_{k,[a^k,n]}(na^ky_k)\frac{f_{k,n}(a^ky_k)}{f_{k,n}(y_k)}\right) \\
      = &\frac{F(any_1,\ldots,a^{l}ny_{l})\left/
            F(ny_1,\ldots,ny_{l})\right.}
            {F(ay_1,\ldots,a^{l}y_{l})\left/
            F(y_1,\ldots,y_{l})\right.} \\
      = &\frac{F([a,n]nay_1,\ldots,[a^{l},n]na^{l}y_{l})}
              {F(nay_1,\ldots,na^{l}y_{l})}
         \frac{ G_n(ay_1,\ldots,a^{l}y_{l})}
              {G_n(y_1,\ldots,y_{l})},
 \end{split}
\end{equation*}
where 
\[
G_n(y_1,\ldots,y_{l})=F(ny_1,\ldots,ny_{l})/F(y_1,\ldots,y_{l}).
\]
By induction using condition (\ref{I:invariance})
\[
\frac{F([a,n]nay_1,\ldots,[a^{l},n]na^{l}y_{l})}
     {F(nay_1,\ldots,na^{l}y_{l})}
\p_{k=1}^{l} \left(\bar{f}_{k,[a^k,n]}(na^{k}y_{k})\right)
\]
is constant $\bt_{l}(\mu_Y)$ a.e.
\end{proof}

\begin{cor}\label{g_n_countable} 
For $k=1,\ldots,l$, $g_{k,n}$ is of type $j-1$, therefore 
the set $\{g_{k,n}\}_{n \in N_i}$  
modulo  $\pi_{j-1}Y$-quasi-coboundaries is countable.
\end{cor}

\begin{proof}
By lemma  \ref{lma:lift}.

\end{proof}

\begin{cor}\label{g_n__r_countable} Suppose for some  $1\le r\le j-1$, $g_{k,n}(y)=g_{k,n}(\pi_ry)$ 
for $k=1,\ldots,l$. Then $g_{k,n}(\pi_ry)$ is of type $r$, and therefore 
the set $\{g_{k,n}(\pi_ry)\}_{n \in N_i}$  
modulo $\pi_{r}Y$-quasi-coboundaries is countable.
\end{cor}

\begin{dsc}{\bf Example.}
Before proceeding with the general proof we describe the proof for 
the case  where $Y$ is a homogeneous space of a $3$-step nilpotent group 
(i.e $N_4=1$, $j-1=3$), $f_k: Y \ra S^1$. Let $n,m \in N_2$. We use the
facts that $[a,m]\in N_3$, therefore 
\[
f_{k,[a,m]}^n f_{k,n} \sim f_{k,n}^{[a,m]}f_{k,[a,m]},\]
and that
for $m_1,m_2 \in N_3$, 
\[f_{k,m_1m_2}\sim f_{k,m_1}^{m_2}f_{k,m_2}.
\]
Recall that by corollary \ref{g_n_countable} there are countably many 
$g_{k,n}$ up to $\pi_3Y$-quasi-coboundaries. Therefore there exists $m \in N_2$
so that for $n$ in a neighborhood of zero in $N_2$
$(g_{k,nm}/g_{k,m})(\pi_3y)$ is a quasi-coboundary:
\[
(g_{k,nm}/g_{k,m})(\pi_3y) \sim L_{k,n,m}(a\pi_3y)/ L_{k,n,m}(\pi_3y)
\]
Fix $n$ and replace  $f_{k,m}(y)$ with 
$f_{k,m}(y) L_{k,n,m}(\pi_3y)$ (this does not effect the commutation relations with $(c,f_{k,c})$ for $c \in N_3$). Now $g_{k,nm}/g_{k,m}$ is a constant.
Then 
\[\begin{split}\frac{f_k(nmy)}{f_k(my)}
&=\frac{f_k(nmy)/f_k(y)}{f_k(my)/f_k(y)}
\sim \frac{f_{k,[a,nm]}(nmay)}{f_{k,[a,m]}(may)}
\frac{f_{k,nm}(am^{-1}my)/f_{k,nm}(m^{-1}my)}{f_{k,m}(am^{-1}my)/f_{k,m}(m^{-1}my)}\\
&\sim f_{k,[n,a]}(many)
\frac{\ti{f}_{k,n}(amy)}{\ti{f}_{k,n}(my)}.
\end{split}\] 
Therefore $g_{k,n}$ may be chosen to be constant (with a proper choice of $f_{k,n}$). Iterating, this holds for all $n \in N_2$. 
To show that we can further modify $f_{k,n}$ so that 
$f_{k,nm} \sim f_{k,n}^mf_{k,n}$ we observe that 
\begin{equation}\label{eq:special_k_nm_2}
\begin{aligned}
  1&=\frac{f_k(nmy)}{f_k(y)}\left/ \frac{f_k(nmy)}{f_k(my)}
  \frac{f_k(my)}{f_k(y)} \right. 
  \sim 
   \frac{f_{k,nm}(ay)\left/ f_{k,n}(may)f_{k,m}(ay)\right.} 
         {f_{k,nm}(y) \left/f_{k,n}(my)f_{k,m}(y)\right.}. \\
\end{aligned}\end{equation}
Now proceed as in lemma \ref{multiplicity_of_f_h}.
As a corollary we get that for $c \in N_3$, $n \in N_2$ we have
$f_{k,n}^cf_{k,c}=f_{k,c}^nf_{k,n}$ 
(same proof as lemma \ref{commuting_center}).
Now let $n \in N_1=N$. We construct $f_{k,n}, g_{k,n}$ as in equation 
(\ref{g_n_still_there}), and as in corollary \ref{g_n_countable} $g_{k,n}$
is of type $3$ and there are 
countably many $g_{k,n}$ up to $\pi_3Y$-quasi-coboundaries.
We first need to establish commutation relations between $(n,f_{k,n})$ and 
$(m,f_{k,m})$
for $m \in N_3$:
\[
 \frac{f_{k,n}(my)}{f_{k,n}(y)} \left/ \frac{f_{k,m}(ny)}{f_{k,m}(y)} \right. 
\]
is a $T$ invariant function and therefore a constant $\gd(k,n,m)$,
which is multiplicative in $m$.
We use this to establish commutation relations between $(n,f_{k,n})$ and 
$(m,f_{k,m})$ 
for $m \in N_2$:
\[\begin{split}
&\frac{T( \frac{f_{k,n}(my)}{f_{k,n}(y)} / \frac{f_{k,m}(ny)}{f_{k,m}(y)}  
      f_{k,[n,m]}(mny))}
     {( \frac{f_{k,n}(my)}{f_{k,n}(y)} / \frac{f_{k,m}(ny)}{f_{k,m}(y)}  
      f_{k,[n,m]}(mny))} \\
=&\gd(k,n,[a,m])\gd(k,a,[m,n])\gd(k,m,[n,a])
  \frac{g_{k,n}(\pi_3 my)}{g_{k,n}(\pi_3y)}\\
=&\gd(k,n,[a,m])\gd(k,a,[m,n])\gd(k,m,[n,a])c(k,n,m)\frac{h_{k,n,m}(\pi_3ay)}{h_{k,n,m}(\pi_3y)},
\end{split}\]
Where $c(k,n,m)$ is multiplicative in $m$ ($g_{k,n}$ is of type $3$ therefore 
satisfies equation \ref{eq:center_N} for $m \in N_2$).
Thus for some eigenfunction $\gp_{k,n,m}(\pi_2y)$
\begin{equation}
      \frac{f_{k,n}(my)}{f_{k,n}(y)}\bar{f}_{k,[n,m]}(nay)
      = \gp_{k,n,m}(y)h_{k,n,m}(\pi_3y)\frac{f_{k,n}(ay)}{f_{k,n}(y)}.
\end{equation}
Let $m,n \in N$. 
As before we find that 
\[
\frac{g_{k,nm}(\pi_3y)}{g_{k,n}(\pi_3my)g_n(\pi_3y)}h_{k,n,[m,a]}(\pi_3 amy)
\sim \frac{T(f_{k,nm}(y)/ f_{k,n}(my)f_{k,m}(y))} 
     {f_{k,nm}(y)/ f_{k,n}(my)f_{k,m}(y)}.
\] 
Same proof as lemma \ref{multiplicity_of_f_h} shows  that 
$f_{k,n},f_{k,m}$ can be chosen so that 
for $m,n$ in a neighborhood of zero in $N$, there exists a function  
$K_{n,m}(\pi_3y) $ such that 
\[
\frac{f_{k,nm}(y)}{f_{k,n}(my)f_{k,m}(y)}=K_{n,m}(\pi_3y).
\]
This implies that $\gd(k,n,c)$ is multiplicative for $n \in N$:
on the one hand
\[
\frac{f_{k,n_1n_2}(cy)}{f_{k,n_1n_2}(y)}
=\gd(k,n_1n_2,c)\frac{f_{k,c}(n_1n_2y)}{f_{k,c}(y)},
\]
while on the other hand
\[\begin{split}
\frac{f_{k,n_1n_2}(cy)}{f_{k,n_1n_2}(y)}
=\frac{f_{k,n_1}(n_2cy)f_{k,n_2}(cy)}{f_{k,n_1}(n_2y)f_{k,n_2}(y)}
=\gd(k,n_1,c)\gd(k,n_2,c)
 \frac{f_{k,c}(n_1n_2y)}{f_{k,c}(y)}.
\end{split}\]
Therefore the constant  
$\gd(k,n,[a,m])\gd(k,a,[m,n])\gd(k,m,[n,a])$ is multiplicative in $m,n$
and must equal $1$
(as $c(k,n,m)$ is locally constant as a function of $n$). 
Thus $c(k,n,m)$ is an eigenvalue of $T$, but it is
multiplicative in $m$ in a neighborhood of zero in $N_2$, therefore  
$c(k,n,m) \equiv 1$. This  implies that $g_{k,n}$ is cohomologous 
to a function on $N/N_2\gC$, and we can chose $f_{k,n}$ so that 
$g_{k,n}$ is lifted from a function on $\pi_2Y$.  
Then $(n,f_{k,n})$ $(m,f_{k,m})$ commute nicely for $m \in N_2$; 
namely $[(n,f_{k,n}),(m,f_{k,m})]=([n,m],cf_{[n,m]})$
for some constant $c$.
Let $n,m \in N$ be so 
$(g_{k,nm}/g_{k,m})(\pi_2y) \sim L_{k,n,m}(a\pi_2y)/ L_{k,n,m}(\pi_2y)$ 
(by corollary \ref{g_n__r_countable} there are countably many 
$g_{k,n}$ up to $\pi_2Y$-quasi-coboundaries) 
Replace  $f_{k,m}(y)$ with 
$f_{k,m}(y) L_{k,n,m}(\pi_2y)$ (this does not effect the commutation relations with $(l,f_{k,l})$ for $l \in N_2$). Now $g_{k,nm}/g_{k,m}$ is a constant.
Computation (using the fact that $(n,f_{k,n}),(nm,f_{k,nm})$ and 
$([a,n^{-1}],f_{k,[a,n^{-1}]})$ commute nicely) shows (as before):
\[
\frac{f_k(nmy)}{f_k(my)}=
\frac{f_k(nmy)/f_k(y)}{f_k(my)f_k(y)}
\sim f_{k,[n,a]}(namy)
\frac{\ti{f}_{k,n}(amy)}{\ti{f}_{k,n}(my)}.
\]
And continue as in lemma \ref{multiplicity_of_f_h} to find that $f_{k,nm}\sim f_{k,n}^mf_{k,m}$.
\end{dsc}

We now proceed with the induction in proposition \ref{P:construction_f_n}. Basically we follow the same procedure. 
\begin{lma} Let $n \in N_i$. The functions $g_n$ and $f_n$ can be chosen so 
that  $g_n$ is a constant function on $Y$.
\end{lma}
\begin{proof}
We already know that for $c \in N_{j-1}$ we have 
$f_{k,n}^cf_{k,c}/ f_{k,c}^nf_{k,n}$ is a constant denoted $\gd(k,n,c)$.
We use induction on $i+1 \le r\le j$ to acquire `good' commuting
relations between  $ (n,f_{k,n}), (m,f_{k,m})$ for $m \in N_r$, 
and to reduce the `level' of $g_{k,n}$. 
Assume that   
\begin{enumerate}
\item $g_{k,n}$ is lifted from $N/N_{r+1}\gC$ and is of type $r+1$. 
\item  For $m \in N_{r+1}$,
\[
f_{k,n}^mf_{k,m}=\gd(k,n,m) f_{k,[m,n]}f_{k,m}^nf_{k,n},
\]
\item\label{eq:star}  For $m \in N_{r+1}$,
\begin{equation*}\begin{split}
\gd(k,n_1n_2,m) f_{k,[m,n_1n_2]}
=&\gd(k,n_1,[n_2,m]) \gd(k,n_1,m) \gd(k,n_2,m) \\
 &f_{k,[[m,n_1],n_2]}^{[n_1,m][n_2,m]} 
              f_{k,[m,n_1]}^{[n_2,m]}f_{k,[m,n_2]},
\end{split}\end{equation*}
and similarly for $\gd(k,n,m_1m_2)$. 
\item $f_{k,n_1n_2}=f_{k,n_1}^{n_2}f_{k,n_2}(y)K_{k,n_1,n_2}(\pi_{r+1}(y))$.
\end{enumerate}
Let $m \in N_r$. Calculations using the induction hypothesis, and commutator
identities lead to 
\begin{equation}\label{eq:Jacobi}\begin{split}
&\frac{T( \frac{f_{k,n}(my)}{f_{k,n}(y)} / \frac{f_{k,m}(ny)}{f_{k,m}(y)}  
      f_{k,[n,m]}(mny))}
     {( \frac{f_{k,n}(my)}{f_{k,n}(y)} / \frac{f_{k,m}(ny)}{f_{k,m}(y)}  
      f_{k,[n,m]}(mny))} \\
=&\gL(k,n,m)
  \frac{g_{k,n}(\pi_{r+1} my)}{g_{k,n}(\pi_{r+1}y)}\\
=&\gL(k,m,n)c(k,n,m)
   \frac{h_{k,n,m}(\pi_{r+1}ay)}{h_{k,n,m}(\pi_{r+1}y)},
\end{split}\end{equation}
where 
\[
\gL(k,m,n_1n_2)= \p_{q=1}^Q \gL(k,p_{1,q}(m,n_1,n_2),p_{2,q}(m,n_1,n_2))
\]
where $p_{1,q}(m,n_1,n_2) \in N_{r}$ is an expression involving $m, m^{-1}$ 
and commutators involving  $m,n_1,n_2,a$ and their inverses, and 
$p_{2,p}(m,n_1,n_2) \in N_i$ is an expression involving $m,n_1,n_2,a$
and there inverses - this is due to condition (\ref{eq:star})
(Same holds for $\gL(k,m_1m_2,n)$),
and $c(k,n,m)$ is multiplicative in $m$ ($g_{k,n}$ is of type $r+1$). As $g_{k,n}$ is countably 
determined up to quasi-coboundaries  $c(k,n,m) \equiv 1$. 
Therefore $g_{k,n}(m\pi_{r+1}y)/g_{k,n}(\pi_{r+1}y)$ is a 
$\pi_{r+1}Y$-coboundary for any
$m \in N_r$, hence $g_{k,n}(\pi_{r+1}y)$ is $\pi_{r+1}Y$-cohomologous to a 
function on $\pi_rY$. 
Therefore we can choose $g_{k,n}, f_{k,n}$ such that $g_{k,n}$ is invariant 
under $m \in N_r$.
This implies that 
\[
f_{k,n}^mf_{k,m}=\gd(k,n,m) f_{k,[m,n]}f_{k,m}^nf_{k,n}
\]
($h_{k,n,m}$ is the constant function $1$). Same calculation as 
in equation (\ref{eq:Jacobi}) gives for any $i<s$, $m \in N_s$   
\[
f_{k,n}^mf_{k,m}=h_{k,n,m}(\pi_{r}y) f_{k,[m,n]}f_{k,m}^nf_{k,n}.
\]
Using this we get, as in lemma \ref{multiplicity_of_f_h}, 
that the functions $f_{k,n}$ can be chosen
so that for $n_1,n_2 \in N_i$ there exists a function 
$K_{k,n_1,n_2} (\pi_ry)$
so that 
\[
f_{k,n_1n_2}(y)=f_{k,n_1}^{n_2}f_{k,n_2}(y)K_{k,n_1,n_2}(\pi_{r}y).
\]
To get the condition on $\gd(k,n,m)$ for $m \in  N_r$: on the one hand
\[
\frac{f_{k,n_1n_2}(my)}{f_{k,n_1n_2}(y)}
=\gd(k,n_1n_2,m)f_{k,[n,m]}(mny)\frac{f_{k,m}(n_1n_2y)}{f_{k,m}(y)},
\]
while on the other hand
\[\begin{split}
\frac{f_{k,n_1n_2}(my)}{f_{k,n_1n_2}(y)}
=&\frac{f_{k,n_1}(n_2my)f_{k,n_2}(my)}{f_{k,n_1}(n_2y)f_{k,n_2}(y)}\\
=&\gd(k,n_1,m)\gd(k,n_2,m)\gd(k,n_1,[n_2,m]) f_{k,[n_1,m]}(mn_1n_2y)\\
 & f_{k,[n_2,m]}(n_1mn_2y)
f_{k,[n_1,[n_2,m]]}([n_2,m]n_1mn_2y)
 \frac{f_{k,m}(n_1n_2y)}{f_{k,m}(y)}.
\end{split}\]
Finally, having good commutation relations between $(n,f_{k,n}),(m,f_{k,m})$
for $n \in N_i$, $m \in N_{i+1}$
(i.e. $[(n,f_{k,n}),(m,f_{k,m})]=([n,m],\gd(k,n,m)f_{k,[n,m]})$),
 we show that for $n \in N_i$, 
the functions $f_{k,n}, g_{k,n}$ can be chosen
so that $g_{k,n}$ is constant. we already know it is lifted from $\pi_{i+1}Y$. 
Let $n,m \in N_i$ be so that  
\[
(g_{k,nm}/g_{k,m})(\pi_{i+1}y) \sim 
h_{k,n,m}(a\pi_{i+1}y)/ h_{k,n,m}(\pi_{i+1}y)
\]
(recall that by corollary \ref{g_n__r_countable} there are countably many 
$g_{k,n}$ up to $\pi_{i+1}Y$-quasi-coboundaries). Replace  $f_{k,m}(y)$ with 
$f_{k,m}(y) h_{k,n,m}(\pi_{i+1}y)$ (this does not affect the commutation relations with $f_{k,p}$ for $p \in N_{i+1}$). Now $g_{k,nm}/g_{k,m}$ is a constant.
Computation shows:
\[\frac{f_k(nmy)}{f_k(my)}=
\frac{f_k(nmy)/f_k(y)}{f_k(my)/f_k(y)}
\sim f_{k,[n,a]}(namy)
\frac{\ti{f}_{k,n}(amy)}{\ti{f}_{k,m}(my)}.
\] 
And continue in lemma \ref{multiplicity_of_f_h} to find that $f_{k,nm}\sim f_{k,n}^mf_{k,m}$.
\end{proof}

\begin{cor}If $i>1$ then for any $c \in N_{j-1}$, $n \in N_i$ 
\[f_{k,n}(cy)f_{k,c}(y)=f_{k,c}(ny)f_{k,n}(y)\]
{\em (i.e. $[(n,f_{k,n}),(c,f_{k,c})]=(1,1)]$)}.
\end{cor}

\begin{proof} The quotient
\[
f_{k,n}(cy)f_{k,c}(y)/f_{k,c}(ny)f_{k,n}(y)
\]
is invariant under translation by $a$ and therefore a constant $\gd(k,n,c)$ 
which  is multiplicative in both coordinates. For $\gc \in  N_{j-1}$,
$f_{\gc}$ is an eigenfunction and therefore invariant under the action
of  $N_i$ for $i>1$. This implies that $\gd(k,n,c\gc)=\gd(k,n,c)$.
Proceed as in lemma 
\ref{commuting_center}.
\end{proof}

\begin{cor}\label{f_n_functional}
If $(n_1,\ldots,n_{l}) \subset N_i^{l}$ 
preserve the ergodic components of $\tau_{l}(Y)$ then 
\begin{equation}\label{eq:invariance}
\frac{F(n_1y_1,\ldots,n_{l}y_{l})}
     {F(y_1,\ldots,y_{l})}
\p_{k=1}^{l}\bar{f}_{k,n_k}(y_{k}) 
\end{equation}
is constant $\bt_{l}(\mu_Y)$ a.e.
\end{cor}

\begin{proof}
Both $([a,n_1],\ldots,[a,n_{l}])$, and  $([a,n_1],\ldots,[a^{l},n_{l}])$
preserve the ergodic components of $\tau_{l}$ (see \ref{dsc:ergodic_com}), and as $g_{k,n}$ is constant
for $n \in N_i$
the function in the left hand of equation (\ref{eq:invariance}) is invariant under $\tau_{l}$, $T_{l}$ 
(this is a calculation, using condition
(\ref{I:commutation_a_n})).
\end{proof}
The proof of proposition \ref{P:construction_f_n} is now complete.
\end{proof}

\begin{pro}\label{pro:construction_group} Let $Y=N/\gC$ be a $(j-1)$-step nilflow, and let 
$(f_1,\ldots,f_l)$ be of type $l$ w.r.t $\bar{\triangle}_l(\mu_Y)$. 
Then for any $k=1,\ldots,l$ the system $Y \times_{f_k} S^1$
can be given the structure of a $j$-step nilflow.
\end{pro}

\begin{proof}
Denote 
\[
\caG_k=\{(n,\gp f_{k,n}): \ n \in N, \gp \ \te{an eigenfunction}, \ f_{k,n} \in \caF_k  \}.
\]
where $\caF_k$ is defined in proposition \ref{P:construction_f_n}.
$\caG_k$ is a group under the multiplication
\[
  (n,f)(m,g)=(nm,f^m g), \qquad  (f^m g)(y)=f(my)g(y).
\]
By proposition \ref{P:construction_f_n}, $\caG_k$ is a $j$-step nilpotent group,
and $(a,f_k) \in \caG_k$. 
Endow $\caG_k$ with the topology:
\[
(n_i,g_i) \ra (n,g) \iff n_i \ra n, \  g_i \converge{L^2(N/\gC)} g.
\]

$\caG_k$ acts transitively and effectively on $X=N/\gC \times S^1$ by:
\[
  (n,f)(y,\gz)=(ny,f(y)\gz).
\]
By a theorem of Montgomery and Zippin 
(see \cite{GOV97} page $88$, Theorem $4.3$) 
it possesses a Lie group structure. This type of construction is carried out 
in \cite{Me90} for the case $j=2$.
\end{proof}

\begin{rmr}
It may also be possible to construct the nilflow is using the constants 
$\gd(k,n,m)$
as was done in \cite{R93},\cite{Le93} for $2$-step nilpotent groups,
and  in \cite{Z02b} for $3$-step nilpotent groups.
\end{rmr}

\begin{lma}\label{lma:f_k_type_j} Let $Y=N/\gC$ be a $(j-1)$-step nilflow, 
and let $f:Y \ra S^1$ be of type $j$. 
Then the system $Y \times_{f} S^1$
can be given the structure of a $j$-step nilflow.
\end{lma}

\begin{proof} By definition, for $k=1,\ldots,l$, there exist integers $m_k$,
with $m_k=1$ for some $k$, and  $(f^{m_1}, \ldots,f^{m_l})$ 
of type $l$ w.r.t $\bt_l(\mu_Y)$. Now use proposition 
\ref{pro:construction_group}.

\end{proof}

A similar proof gives:
\begin{lma}\label{lma:f_g_type_j} Let $Y=N/\gC$ be a $(j-1)$-step nilflow, 
and let $f,g:Y \ra S^1$ be of type $j$. 
Then the system $Y \times_{fg} S^1$
can be given the structure of a $j$-step nilflow.
\end{lma}

\begin{cor}\label{cor:pronillow} Let $Y$ be a $(j-1)$-step pro-nilflow, 
and let $(f_1,\ldots,f_l)$ be of type $l$ w.r.t $\bar{\triangle}_l(\mu_Y)$. 
Then for any $k=1,\ldots,l$ the system $Y \times_{f_k} S^1$
can be given the structure of a $j$-step pro-nilflow.
\end{cor}
\begin{proof} By lemma \ref{finite_dim_nil} and corollary \ref{cor:any_l},
$f_k$ is cohomologous to a cocycle $\ti{f}_k$
lifted from a $(j-1)$-step nilflow $(N/\gC,a)$. Furthermore, 
there exist $(g_1,\ldots,g_l)$
with $g_k$ of type $j-1$, 
such that $(\ti{f}_1g_1,\ldots,\ti{f}_lg_l)$ 
is of type $l$ w.r.t  $\bar{\triangle}_l(\mu_{N/\gC})$. By proposition \ref{pro:construction_group}, $N/\gC \times_{\ti{f}_kg_k} S^1$
can be given the structure of a $j$-step nilflow. By
lemma \ref{lma:f_k_type_j}, $N/\gC \times_{g_k^{-1}} S^1$ 
can be given the structure of a $j$-step nilflow.
By the construction in proposition \ref{pro:construction_group},
$N/\gC \times_{\ti{f}_kg_kg_k^{-1}} S^1$
can be given the structure of a $j$-step nilflow.
\end{proof}

\begin{dsc}\label{proof:type_j_pronil}{\em proof of \ref{thm:main}(\ref{item:type_j_pronil})}

Let  $X=Y_{j}(X)\times_{\gr} H$, where $H$ is a compact abelian group,
and for any $\gx \in \hat{H}$, there exists 
$(\gx_1,\ldots,\gx_l) \in \hat{H}^l$, with $\gx=\gx_k$ for some $k$,
and $(\gx_1 \circ \gr,\ldots,\gx_l \circ \gr)$ of type $\vec{a}$ w.r.t
$\bt_{\vec{a}}(Y_j(X))$. By corollary \ref{cor:pronillow} the system
 $X=Y_{j}(X)\times_{\gx \circ \gr} S^1$ is isomorphic to a $j$-step
pro-nilflow.
By Pontryagin duality, $H \hookrightarrow (S^1)^{\hat{H}}$. 
The system 
$Y \times_{\gr} H$ is therefore a `join' of factors of the form 
$Y \times_{\gx \circ \gr} S^1$ where $\gx$ ranges over $\hat{H}$.
\end{dsc}

\begin{lma}\label{gl_trivial}  Let $Y$ be a $(j-1)$-step pro-nilflow, and let 
$(f_1,\ldots,f_l)$ be of type $l$ w.r.t $\bar{\triangle}_l(\mu_Y)$. 
If for some $k$, $\gl_{k,u} \equiv 1$ (see lemma \ref{lma:g_u_k_l}) in 
a neighborhood of zero in $H_{j}$, then the system $Y \times_{f_k} S^1$
can be given the structure of a $(j-1)$-step pro-nilflow.
\end{lma}

\begin{proof}
By the induction hypothesis \ref{thm:main}(\ref{item:abelian_ext_type_j}),
$H_{j}$ is connected ($j \ge 1$)
therefore we may choose $\gl_{k,u} \equiv 1$
on $ H_{j}$. By corollary \ref{hom}, $f_k$ is cohomologous to a cocycle 
$\ti{f_k}$ lifted from   $Y^0_{j-1}$. By lemma \ref{lma:lift},
$\ti{f_k}$ is of type $j-1$. By 
the induction hypothesis \ref{thm:main}(\ref{item:type_j_pronil}),
$Y_{j-1} \times_{\ti{f_k}} S^1$ can be given the structure of a 
$(j-1)$-step pro-nilflow. 
\end{proof}

\begin{lma}\label{H_j_connected} Let  $X=Y_{j}(X)\times_{\gr} H$, where $H$ is a compact abelian group, and $\gr$ of type $j$. Then $H$ 
is connected.
\end{lma}
\begin{proof} 
We show that $\hat{H}$ has no elements of finite order.
Suppose for some $\gx \in \hat{H}$, and some 
$l>0$, $\gx^l=1$. $\gx \circ \gr$ satisfies equation 
(\ref{center_eq}). By lemma 
\ref{lma:mult_of_f_u} the function $\gl_{k,u}$ 
is multiplicative in a 
neighborhood of zero in $H_{j}$. $\gl_{k,u}^l$ is an eigenvalue,
and as $H_{j}$ is connected $\gl_{k,u} \equiv 1$ in a neighborhood of zero. 
By lemma \ref{gl_trivial} the system 
$Y_j'=Y_j(X) \times_{\gx \circ \gr} S^1$, 
which is a factor of $X=Y_j(X) \times_{\gr} H$, can be given the structure
of a $(j-1)$-step pro-nilflow. By corollary
\ref{cor:ucf_pronil},
$Y_j(Y'_j)=Y'_j$,  in contradiction to  
$Y_j(X)$ being the $j$-u.c.f of $X$.
\end{proof}

\begin{lma}\label{lma:abelian}
Let $X$ be a group extension of $Y_{j}(X)$;
i.e $X={Y}_{j}(X) \times_{\gs} G$. Then  
$Y_{j+1}(X)$ is an abelian extension of ${Y}_{j}(X)$ by a cocycle of type $j$,
and therefore can be given the structure of a $j$-step pro-nilflow.
\end{lma}

\begin{proof}
The proof is a straightforward generalization of lemmas  $9.1,9.2$ in 
\cite{FuW96} 
(this is done for the case $j=3$ in \cite{Z02b}). 
We outline the steps. Any ergodic component of 
$\bar{\triangle}_{\vec{a}}(\mu_{X})$ projects
onto an ergodic component of  $\bar{\triangle}_{\vec{a}}(\mu_{Y_{j}(X)})$. 
The fact that 
$\tau_{\vec{a}}(T)$, and $T_{j+1}(T)$ commute implies that the 
Mackey groups associated with 
different ergodic component of $\bar{\triangle}_{\vec{a}}(\mu_{Y_{j}(X)})$ 
are conjugate for a.e. ergodic component (lemma \ref{lma:conjugate}). 
Denote $[M_{\vec{a}}]$ the conjugacy class, where 
the  group $M_{\vec{a}} \subset G^{j+1}$. 
One then uses the fact that the projection
of $M_{\vec{a}}$ on any $j$ coordinates is full (i.e. $G^{j}$)
to show that 
$[G,G]^{j+1} \subset M_{\vec{a}}$. More specifically one shows that
there exists an abelian group $K_{\vec{a}}$ and homomorphisms 
$\gp_{\vec{a},i}:G \ra K_{\vec{a}}$ so that
\[
 M_{\vec{a}}=\{(g_1,\ldots,g_{j+1}|\gp_{\vec{a},i}(g_1)\ldots \gp_{\vec{a},i}(g_{j+1})=1\}.
\]

We return to the average in (\ref{D}). By \ref{P:reduction_abelian} we can 
replace 
\[
f_1 \otimes \ldots \otimes f_{j+1}(y_1,g_1,\ldots,y_{j+1}g_{j+1})
\]
by 
\[
\int f_1 \otimes \ldots \otimes f_{j+1}
(y_1,g_1m_1,\ldots,y_{j+1},g_{j+1}m_{j+1})
dm_{M_{\vec{a}}}(m_1,\ldots,m_{j+1})
\]
where $dm_{M_{\vec{a}}}$ is the Haar measure on the Mackey group $M_{\vec{a}}$.
As $[G,G]^{j+1} \subset M_{\vec{a}}$ we can replace $f_k$, 
for $k=1,\ldots,j+1$,
by $\int f_1(y,gg')dm_{[G,G]}(g')$. Thus  
${Y}_{j}(X) \times_{\gr} G/[G,G]$ is characteristic for the scheme $\vec{a}$,
for any $\vec{a}$.

  Let   $K_0=\cap_{k,\vec{a}}$ ker$\gp_{\vec{a},k}= \{1\}$.
Let $\ti{G}=G/K_0$, and  let $H=\ti{G}/[\ti{G},\ti{G}]$. Then similarly
${Y}_{j}(X) \times_{\gr} H$ is characteristic for the scheme $\vec{a}$,
for any $\vec{a}$. We will show that $\gr$ is of type $j$. Then by 
\ref{proof:type_j_pronil} it can be given the structure of a $j$-step pro-nilflow, and 
by corollary \ref{cor:ucf_pronil} it the $j+1$ universal characteristic factor.

Denote $Y=Y_j(X)$.
Then by equation (\ref{eq:conditional}) 
\[
\bar{\triangle}_{\vec{a}}(\mu_{X})=\bar{\triangle}_{\vec{a}}(\mu_{Y})
 \times m_{H}^{j+1}.
\] 
For each ergodic component of 
$\bar{\triangle}_{\vec{a}}(\mu_{Y})$ the ergodic components
of $\bar{\triangle}_{\vec{a}}(\mu_X)$ are given by the Mackey group 
$M_{\vec{a}} \subset H^{j+1}$. 
Above a.e. ergodic component $W_{\vec{a},y}$ (by the discussion in \ref{dsc:ergodic_com}
the ergodic components of $\bar{\triangle}_{\vec{a}}(\mu_{Y})$ are 
parametrized by $Y$) we have a $H^{j+1}$- extension
by the cocycle 
\[
\ti{\gr}_{\vec{a}}=(\gr^{(a_1)}(y_1),\gr^{(a_2)}(y_2),\ldots,\gr^{(a_{j+1})}(y_{j+1}))
: W_{\vec{a},y} \ra H^{j+1}.
\] 
By Theorem \ref{mackey_thm}
there exists a function $\gf:W_{\vec{a},y} \ra H^{j+1}$
such that
\[
\gf_{\vec{a}}(\tau_{\vec{a}}(y_1,\ldots,y_{j+1}))\ti{\gr}_{\vec{a}}(y_1,\ldots,y_{j+1})
\gf_{\vec{a}}^{-1}(y_1,\ldots,y_{j+1}) \in M_{\vec{a}}
\] 
Applying the foregoing characterization of $M_{\vec{a}}$, there exists 
an abelian group $K_{\vec{a}}$ and homomorphisms $\gp_{\vec{a},i}:H \ra K_{\vec{a}}$
\begin{equation}\label{eq:mackey_condition}
 \p_{k=1}^{j+1} \gp_{\vec{a},k}\circ \gr^{(a_k)} (y_k)= F_{\vec{a}}(\tau_{\vec{a}}(y_1,\ldots,y_{j+1}))
F_{\vec{a}}^{-1}(y_1,\ldots,y_{j+1})
\end{equation}
Where 
\[
F_{\vec{a}}(y_1,\ldots,y_{j+1})=\p_{k=1}^{j+1} \gp_{\vec{a},i} \circ 
\gf_{\vec{a},i} (y_1,\ldots,y_{j+1}) \in K_{\vec{a}}.
\]
Let $\gx \in  \hat {K}_{\vec{a}}$. Applying $\gx$
to equation (\ref{eq:mackey_condition}) we get 
\begin{equation}\label{eq:mackey}
  \p_{k=1}^{j+1} \chi \circ \gp_{\vec{a},k} \circ (\gr)^{(a_k)} (y_k)
     =\frac{\tau_{\vec{a}} F_{\vec{a},y,\chi}(y_1,\ldots,y_{j+1})}
            {F_{\vec{a},y,\chi}(y_1,\ldots,y_{j+1})}.
\end{equation}
Where $F_{\vec{a},y,\gx}:W_{\vec{a},y} \ra S^1$. By ergodicity of 
$\tau_{\vec{a}}$ on $W_{\vec{a},y}$, $F_{\vec{a},y,\chi}$
is unique up to a constant multiple. 
By proposition \ref{g_x_constant} there is a measurable choice of $F_{\vec{a},y,\gx}$, 
so that  equation  ({\ref{eq:mackey}})
holds  $\bar{\triangle}_{\vec{a}}(\mu_{Y})$ a.e. 
Finally, as $\cap_{k=1}^{j+1}$ ker $\gp_{\vec{a},k} = \{1\}$,
the characters $\gx \circ \gp_{\vec{a},k}$ where $k=1,\ldots,j+1$, $\vec{a} \in \BZ^{j+1}$, and $\gx\in \hat{K}_{\vec{a}}$
span $\hat{H}$.
\end{proof}

\begin{dsc}\label{cor:presentation_abelian} {\em proof of theorem \ref{thm:main}(\ref{item:abelian_ext_type_j})}. 

If $X$ is a $j$-step pro-nilflow then $X$ is an abelian extension of $Y_j(X)$.
By corollary \ref{cor:ucf_pronil}, $X=Y_{j+1}(X)$. By 
lemma \ref{lma:abelian} it is an extension of $Y_j(X)$ by a cocycle of 
type $j$.
\end{dsc}

\begin{cor}\label{cor:finite_pronil_is_nil}
Any $j$-step pro-nilflow $Y$ can be presented as a tower of abelian extensions
$H_1 \times_{\gs_1} H_2 \times \ldots \times_{\gs_j} H_{j+1}$ where
$\gs_k$ for $k=1,\ldots,j$ is of type $k$. If in this presentation $H_1$
is a cyclic group, 
and for $k>1$,  $H_k$ is a finite dimensional torus, then $Y$ is a nilflow.
\end{cor}

\begin{proof} The first part is clear. The second part follows from 
the construction in proposition
\ref{pro:construction_group}.
\end{proof}

\begin{dsc}\label{proof:factor_nil} {\em proof of \ref{thm:main}(\ref{item:factor_nil})}.

Let 
\[{X}=(N/\gC,T)=(\lim_{\leftarrow} N_i/\gC_i,a_i)\]
be a $j$-step pro-nilflow, and let ${W}$
be a factor. Let ${Y}_j=\lim_{\leftarrow} M_i/\gL_i$ be the $j$-u.c.f 
of ${X}$. Then ${X}={Y}_j \times_{\gs_j} H$
where $H$ is a compact abelian group,  and $\gs_j$ is of type $j$.
Let $K$ be a compact abelian group of m.p.ts acting on $X^0$ and commuting 
with the action of $T$. We show that $K$ commutes with the action of $N$.
By corollary \ref{natural}, any $k \in K$ induces  a map from $Y_j$ to itself,
also denoted $k$ by abuse of notation.
The action of $K$ is given by $k(y,h)=(ky,\gr_k(y,h))$. 
We first show that
$k$ preserves the skew product structure.
As $k,T$ commute:
\[
\gr_k(T(y,h))=\gr_k(y,h)+\gs_j(ky).
\]    
Denote 
\[ F_k(y,h)=\gr_k(y,h)-h.\]
Then
\[
TF_k(y,h)-F_k(y,h)=\gs_j(ky)-\gs_j(y)
\]
Let $\gx$ be a character of $H$, then 
 \[
\frac{\gx \circ \gs_j(ky)}{\gx \circ\gs_j(y)}
=\frac{T\gx \circ F_{k}(y,h)}{\gx \circ F_{k}(y,h)}.
\]
As $\gs_j$ is of type $j$, $\gx \circ \gs_j$ is lifted from $M_i/\gL_i$ 
for some $i$. Let $p:Y_j \ra M_i/\gL_i$ be the projection.
By induction $k_j$ commutes with the action of $M= \lim_{\leftarrow} M_i$.
By proposition \ref{P:construction_f_n}, as $k$ commutes with the action of 
$M= \lim_{\leftarrow} M_i$, 
there exist functions measurable $f_{\gx,k}:M_i/\gL_i \ra S^1$ such that 
 \[
\frac{\gx \circ \gs_j( ky)}{\gx \circ\gs_j( y)}
=\gl_{k}\frac{Tf_{\gx,k}(p y)}{f_{\gx,k}(p y)}.
\]
Therefore
\[
\frac{\bar{f}_{\gx,k}(Tp y)\gx \circ F_k(T(y,h))}
{\bar{f}_{\gx,k}(p y)\gx \circ F_k(y,h)}
=\gl_{k}
\]
 Therefore $f_{\gx,k}(p y)\gx \circ F_k(y,h)$ is an eigenfunction 
$\gp(\pi_2y)$ (it is defined on the Kronecker factor). As $j\ge3$,
$F_k(y,h)$ depends only on $y$. Therefore 
\[
\gx \circ \gr_k(y,h)=\gx(h)\gp(\pi_2y)f_{\gx,k}(p y).
\]
This implies that the action of $k$ induces an action on $N_i/\gC_i$
for all $i$, that commutes with the action of $T$, and by \cite{P73}
Theorem $4.3$
it commutes with the action of $N_i$ (the proof in \cite {P73} is 
for $(N/\gC,a)$ where $N$ is connected,
but same proof holds in the case where  $N$ is generated by $a$ and the 
connected component of the identity).
We obtain the result inductively, using the fact that 
${X}$ has generalized discrete spectrum mod $\caD$ of
finite type (see \cite{P73}), and is therefore obtained from ${W}$
by a finite series of abelian extensions. 

\end{dsc}

\begin{dsc}\label{proof:main} {\em proof of theorem \ref{thm:main}(\ref{item:main})}.

By \ref{proof:iso}, 
${Y}_{j+1}(X)$ is an isometric extension of the factor ${Y}_{j}(X)$. 
By  the discussion in \ref{isometric_extensions}, 
${Y}_{j+1}(X)$ is  of the form ${Y}_{j}(X) \times_{\gs} G/L$, 
where $G/L$ is a homogeneous
space of a compact metric group $G$. By lemma \ref{isometric_to_group} 
we may assume that the extension $X'={Y}_{j}(X) \times_{\gs} G$ is an 
ergodic group extension. As ${Y}_{j}(X)$ is a factor of ${Y}_{j}(X')$, by 
lemma \ref{lma:intermediate}, $X'$ is  group extension of ${Y}_{j}(X')$; i.e.
$X'={Y}_{j}(X') \times_{\gs'} G'$. By corollary \ref{natural}, 
the factor map $X' \ra X$  induces a map 
between their $j+1$-u.c.fs.  By \ref{proof:factor_nil},  it is enough to
show that ${Y}_{j+1}(X')$ has the structure
f a $j$-step pro-nilflow. By lemma \ref{lma:abelian} we are done.
\end{dsc}

\begin{dsc}\label{proof:set_type_j_to_type_j} {\em proof of theorem 
             \ref{thm:main}(\ref{item:set_type_j_to_type_j})}.

By  \ref{cor:pronillow}, 
$X=Y \times_{f_k} S^1$ can be given the structure of a $j$-step 
pro-nilflow, thus $Y_{j+1}(X)=X$. If 
$Y_{j}(X)=Y$, then by \ref{cor:presentation_abelian}, 
$f_k$ is of type $j$. Otherwise,  $Y$ is a proper factor of  $Y_{j}(X)$ and 
therefore $Y_{j+1}(X)$ is an extension
of $Y_{j}(X)$ by a proper closed subgroup $G$ of $S^1$. 
By \ref{cor:presentation_abelian} and \ref{H_j_connected}, $G$ must be trivial,
thus $Y_{j+1}(X)=Y_{j}(X)$, which is a $(j-1)$-step pro-nilflow. 
This implies that we can chose $f_{u,k},\gl_{u,k}$ in equation 
(\ref{center_eq}), with  $\gl_{u,k} \equiv 1$ 
(otherwise by proposition \ref{P:construction_f_n} and the construction
in \ref{pro:construction_group}, 
we increase the level of nilpotency). By corollary \ref{hom}, $f_k$
is cohomologous to a function $f'_k$ on $Y_{j-1}(Y)$. The system 
$Y_{j-1} \times_{f'_k} S^1$ is a factor of $X$ and therefore a $j-1$-step
pro-nilflow. By the induction hypothesis 
\ref{main}(\ref{item:set_type_j_to_type_j}) $f'_k$ is of type
$j-1$, therefore $f_k$ is of type $j$.  
\end{dsc}

\begin{dsc} {\em proof of theorem \ref{thm:main}(\ref{item:type_j_to_type_j})}.

If $f:Y \ra H$ is of type $j$ then  for any $\gx \in \hat{H}$ there exists
$(\gx_1, \ldots ,\gx_l) \in \hat{H}^l$ with $\gx=\gx_k$ for some $k$, and
 $(\gx_1\circ f, \ldots ,\gx_l\circ f)$ is of type $\vec{a}$ w.r.t 
$\bt_{\vec{a}}(Y)$. By \ref{proof:set_type_j_to_type_j}, $\gx_l\circ f$
is of type $j$.
\end{dsc}

\begin{dsc} {\em proof of theorem \ref{thm:main}(\ref{item:group})}.

By  lemma \ref{lma:f_g_type_j}, $Y \times_{fg} S^1$ can be given the 
structure of a  pro-nilflow (it is clear from the proof of 
\ref{finite_dim_nil} that we can have the functions 
$f,g$ lifted form the same $j-1$-step nilflow). 
As in \ref{proof:set_type_j_to_type_j}, $fg$ is of type $j$.
\end{dsc}

\begin{dsc}  {\em proof of theorem \ref{thm:main}(\ref{eq:convergence})}. 

This follows from \ref{proof:main}  and \ref{dsc:ergodic_com}.
\end{dsc}

\end{proof}

\end{document}